%% file: dsmanifolds-revised-2.tex
\documentclass[11pt,a4paper]{article}
\usepackage{graphicx}
\usepackage{amsmath}
\usepackage{amsfonts}
\usepackage{dsfont} 
\usepackage{amssymb}
\usepackage{stmaryrd} 
\usepackage[normalem]{ulem}
\usepackage{enumerate}
\usepackage[hidelinks]{hyperref}
\usepackage{lscape}
\usepackage{longtable}
\usepackage{mathrsfs}
\usepackage{diagbox} 
\usepackage{rotating}
\usepackage{booktabs} 
\usepackage{tabularx}
\usepackage{subcaption}
\usepackage{array}
\usepackage{url}
\usepackage{comment}
\usepackage{cite} 
\usepackage[dvipsnames]{xcolor}
\usepackage{colortbl} 
\usepackage{arydshln} 
\usepackage{yfonts}
\usepackage{tikz}
\usetikzlibrary{calc}

\usepackage{todonotes}
\usepackage[labelfont=sl,textfont=sl]{subcaption} 
\usepackage[font={small,sl},labelsep=period]{caption} 
\usepackage[ruled,vlined,noline,linesnumbered]{algorithm2e}

\newtheorem{theorem}{Theorem}[section]

\newtheorem{corollary}{Corollary}[section]

\newtheorem{definition}{Definition}[section]

\newtheorem{example}{Example}[section]

\newtheorem{lemma}{Lemma}[section]

\newtheorem{proposition}{Proposition}[section]
\newtheorem{remark}{Remark}[section]

\newtheorem{assumption}{Assumption}[section]

\newenvironment{proof}[1][Proof]{\textbf{#1.} }{\hfill$\Box$\vspace{\baselineskip}}
\numberwithin{equation}{section}
\newcolumntype{M}[1]{>{\centering\arraybackslash}m{#1}}

\newcommand{\R}{\mathbb{R}}

\newcommand{\calL}{\mathcal{L}}
\newcommand{\calS}{\mathcal{S}}
\newcommand{\calU}{\mathcal{U}}

\DeclareMathOperator{\spann}{span}
\DeclareMathOperator{\proj}{proj}

\DeclareMathOperator{\Skew}{Skew}

\newcommand{\T}{\mathrm{T}}

\newcommand{\norm}[1]{\left\|#1\right\|}

\newcommand{\inner}[2]{\left\langle{#1},{#2}\right\rangle}

\newcommand{\SO}{\mathrm{SO}}

\newcommand{\bbS}{\mathbb{S}}

\newcommand{\Retr}{\mathrm{R}}
\newcommand{\grad}{\mathrm{grad}}

\newcommand{\calM}[0]{\mathcal{M}}
\newcommand{\M}[0]{\mathcal{M}}

\newcommand{\Mlinear}[0]{F}

\newcommand{\Msphere}[0]{\mathcal{M}_s}
\newcommand{\Mrotsync}[0]{\mathcal{M}_{rs}}

\newcommand{\cm}{\mathrm{cm}}
\newcommand{\D}{D} 
\newcommand{\DM}{\mathcal{D}} 
\newcommand{\B}{B} 
\newcommand{\BM}{\mathcal{B}} 
\newcommand{\fP}{P} 
\newcommand{\PM}{\mathcal{P}} 

\newcommand{\calQ}{\mathcal{Q}}
\newcommand{\calE}{\mathcal{E}}

\newcommand{\cardpss}{r}

\DeclareRobustCommand{\pssa}{\D^{\mathrel{\text{\protect\begin{tikzpicture}[scale=0.1, baseline=-0.5ex]
			\draw[-] (0,0) -- (0,1);
			\draw[-] (0,0) -- (0,-1);
			\draw[-] (0,0) -- (1,0);
			\draw[-] (0,0) -- (-1,0);
			\draw[] (0,0) circle (1.0);
		\end{tikzpicture}}}}}
\DeclareRobustCommand{\pssb}{\D^{\mathrel{\text{\protect\begin{tikzpicture}[scale=0.1, baseline=-0.5ex]
		\draw[-] (0,0) -- (0,1);
		\draw[-] (0,0) -- (1,0);
		\draw[-] (0,0) -- (-0.707,-0.707);
		\draw[] (0,0) circle (1.0);
\end{tikzpicture}}}}}
\DeclareRobustCommand{\pssc}{\D^{\mathrel{\text{\protect\begin{tikzpicture}[scale=0.1, baseline=-0.5ex]
	\begin{scope}[rotate = 30]
			\draw[-] (0,0) -- (1,0);
			\draw[-] (0,0) -- (-1/2,0.866);
			\draw[-] (0,0) -- (-1/2,-0.866);
	\end{scope}
	\draw[] (0,0) circle (1.0);
\end{tikzpicture}}}}}

\newcommand{\pssaintr}[1]{%
	\D^{\mathrel{\text{\protect\begin{tikzpicture}[scale=0.1, baseline=-0.5ex]
			\draw[-] (0,0) -- (0,1);
			\draw[-] (0,0) -- (0,-1);
			\draw[-] (0,0) -- (1,0);
			\draw[-] (0,0) -- (-1,0);
			\draw[] (0,0) circle (1.0);
		\end{tikzpicture}}}} (#1)}

\newcommand{\pssbintr}[1]{%
	\D^{\mathrel{\text{\protect\begin{tikzpicture}[scale=0.1, baseline=-0.5ex]
		\draw[-] (0,0) -- (0,1);
		\draw[-] (0,0) -- (1,0);
		\draw[-] (0,0) -- (-0.707,-0.707);
		\draw[] (0,0) circle (1.0);
\end{tikzpicture}}}} (#1)}

\newcommand{\psscintr}[1]{%
	\D^{\mathrel{\text{\protect\begin{tikzpicture}[scale=0.1, baseline=-0.5ex]
	\begin{scope}[rotate = 30]
			\draw[-] (0,0) -- (1,0);
			\draw[-] (0,0) -- (-1/2,0.866);
			\draw[-] (0,0) -- (-1/2,-0.866);
	\end{scope}
	\draw[] (0,0) circle (1.0);
\end{tikzpicture}}}} (#1)}

\DeclareRobustCommand{\pssaproj}[2]{%
\if\relax\detokenize{#1}\relax \PM \else\fP_{#1}\fi(\D^{\mathrel{\text{\protect\begin{tikzpicture}[scale=0.1, baseline=-0.5ex]
			\draw[-] (0,0) -- (0,1);
			\draw[-] (0,0) -- (0,-1);
			\draw[-] (0,0) -- (1,0);
			\draw[-] (0,0) -- (-1,0);
			\draw[] (0,0) circle (1.0);
		\end{tikzpicture}}}}\if\relax\detokenize{#2}\relax \else (#2)\fi )  }%

\newcommand{\pssbproj}[2]{%
\if\relax\detokenize{#1}\relax \PM \else\fP_{#1}\fi(\D^{\mathrel{\text{\protect\begin{tikzpicture}[scale=0.1, baseline=-0.5ex]
		\draw[-] (0,0) -- (0,1);
		\draw[-] (0,0) -- (1,0);
		\draw[-] (0,0) -- (-0.707,-0.707);
		\draw[] (0,0) circle (1.0);
\end{tikzpicture}}}}\if\relax\detokenize{#2}\relax \else (#2)\fi )  }%

\newcommand{\psscproj}[2]{%
\if\relax\detokenize{#1}\relax \PM \else\fP_{#1}\fi(\D^{\mathrel{\text{\protect\begin{tikzpicture}[scale=0.1, baseline=-0.5ex]
	\begin{scope}[rotate = 30]
			\draw[-] (0,0) -- (1,0);
			\draw[-] (0,0) -- (-1/2,0.866);
			\draw[-] (0,0) -- (-1/2,-0.866);
	\end{scope}
	\draw[] (0,0) circle (1.0);
\end{tikzpicture}}}}\if\relax\detokenize{#2}\relax \else (#2)\fi )  }%

\newcommand{\eps}{\epsilon}
\newcommand{\alphamax}[0]{\alpha_{\mathrm{max}}}
\newcommand{\gammainc}[0]{\gamma_{\mathrm{inc}}}
\newcommand{\gammadec}[0]{\gamma_{\mathrm{dec}}}
\newcommand{\flow}[0]{f_{\mathrm{low}}}
\newcommand{\Dmax}[0]{d_{\mathrm{max}}}
\newcommand{\Dmin}[0]{d_{\mathrm{min}}}
\newcommand{\Skj}{\calS_{j_1}(k_0)}
\newcommand{\Ukj}{\calU_{j_1}(k_0)}
\newcommand{\cmx}{\mathrm{cm}_{x}}
\newcommand{\isomintr}[1]{\varphi \if\relax\detokenize{#1}\relax\else {_{#1}}\fi}

\newcommand{\rev}[1]{#1}
\newcommand{\rrev}[1]{{{\color{black}{#1}}}} 
\newcommand{\prev}[1]{{{\color{black}{#1}}}} 

\begin{document}

\title{Complexity guarantees and polling strategies for Riemannian direct-search methods}
\author{
	Bastien Cavarretta\thanks{LAMSADE, CNRS, Universit\'e Paris Dauphine-PSL,
		Place du Mar\'echal de Lattre de Tassigny, 75016 Paris, France
		(\texttt{bastien.cavarretta@dauphine.psl.eu}).}
	\and
	Florentin Goyens\thanks{ICTEAM Institute, UCLouvain, Louvain-la-Neuve, Belgium. ORCID 0000-0001-5531-5673 (\texttt{florentin.goyens@uclouvain.be}).}
	\and
	Cl\'ement W. Royer\thanks{LAMSADE, CNRS, Universit\'e Paris Dauphine-PSL,
		Place du Mar\'echal de Lattre de Tassigny, 75016 Paris, France.
		ORCID 0000-0003-2452-2172
		(\texttt{clement.royer@lamsade.dauphine.fr}).}
	\and
	Florian Yger\thanks{LITIS, INSA Rouen Normandie, 685 avenue de l'Université, 76801 Saint-Etienne-du-Rouvray (\texttt{florian.yger@insa-rouen.fr}).}
}
\date{}
\maketitle
\thispagestyle{empty}
\footnotesep=0.4cm

\begin{abstract}
	Direct-search algorithms are derivative-free optimization techniques that operate by polling the variable space along specific directions forming
	positive spanning sets (PSSs). When the problem variables are constrained
	to lie on a Riemannian manifold, polling must be performed along tangent
	directions. Although Riemannian variants of direct search have already
	been proposed and endowed with asymptotic guarantees, a proper
	generalization of PSSs on manifolds remains to be investigated. In
	particular, a measure of quality for those PSSs is required to obtain
	complexity bounds for direct search.

	In this paper, we derive complexity guarantees for a class of Riemannian
	direct-search techniques, and study two ways of generating positive
	spanning sets in tangent spaces. We pay particular attention \rev{to} the
	unit hypersphere case, for which we establish that genera\-ting directions
	directly within the tangent space leads to better complexity properties
	than projecting PSSs from the ambient space onto the tangent space. Our
	numerical experiments highlight the impact of dimension and codimension
	in more general settings.
\end{abstract}

\section{Introduction}
\label{sec:intro}

Direct-search methods form one of the main classes of derivative-free
optimization algorithms. Despite their apparent simplicity, these methods
enjoy solid theoretical foundations that have been developed over several
decades~\cite{CAudet_2014,CAudet_WHare_2017,ARConn_KScheinberg_LNVicente_2009,
	KJDzahini_FRinaldi_CWRoyer_DZeffiro_2025,TGKolda_RMLewis_VTorczon_2003}. At
every iteration of a direct-search algorithm, the variable space is explored
\rev{through} a set of polling directions. \rev{In the unconstrained setting,
	these directions are typically required to form a positive spanning
	set~\cite{CDavis_1954}, meaning that they generate the variable space via
	nonnegative linear combinations. Thanks to this property, a direct-search
	method can be endowed with convergence~\cite{TGKolda_RMLewis_VTorczon_2003}
	and worst-case complexity guarantees~\cite{LNVicente_2013}. The latter depend}
upon the number of directions used per iteration and the quality of those
directions, quantified by the cosine measure of the polling
set~\cite{TGKolda_RMLewis_VTorczon_2003}. Computing cosine
measures~\cite{RGRegis_2016,RGRegis_2021,WHare_GJarryBolduc_2020,
	WHare_SSun_2025} and identifying positive spanning sets with good cosine
measures~\cite{WHare_GJarryBolduc_CPlaniden_2023,
	WHare_GJarryBolduc_SKerleau_CWRoyer_2024} have thus become growing areas of
interest in the direct-search community.

In presence of constraints, the polling directions must conform to the
geometry of the feasible set around the current iterate. This is particularly
critical for nonrelaxable constraints, that must be satisfied in order for the
optimization objective to be evaluated~\cite{SLeDigabel_SMWild_2024}. A
typical example of such constraints are bounds and linear constraints, for
which a number of direct-search variants have been
proposed~\cite[Section 5.1.1]{KJDzahini_FRinaldi_CWRoyer_DZeffiro_2025}. In
particular, complexity guarantees can be derived provided the polling directions
are positive spanning sets of certain cones defined by nearby
constraints~\cite{SGratton_CWRoyer_LNVicente_ZZhang_2019}. The special case
of linear equality constraints leads to considering directions within a
linear subspace, for which the notions of positive spanning sets and cosine
measure can be formalized and connected to complexity
bounds~\cite{CAudet_WHare_GJarryBolduc_2025,
	SGratton_CWRoyer_LNVicente_ZZhang_2019,LLiu_XZhang_2006}.

Meanwhile, fueled by the development of manifold
optimization~\cite{PAAbsil_RMahony_RSepulchre_2008,NBoumal_2023}, several
direct-search algorithms have been proposed to tackle problems where the
constraints define a Riemannian manifold~\cite{DWDreisigmeyer_2006,
	DWDreisigmeyer_2007a,VKungurtsev_FRinaldi_DZeffiro_2024}. In that setting,
polling directions must be generated within the tangent space corresponding
to the current iterate, which potentially changes as the algorithm proceeds.
\rev{A possible strategy recently employed by Kungurtsev et	al.~\cite{VKungurtsev_FRinaldi_DZeffiro_2024} consists} in projecting
positive spanning sets from the ambient Euclidean space onto the Riemannian
manifold. Provided the resulting directions are of sufficiently good quality
within the tangent spaces, convergence of the corresponding method can be
established. Nevertheless, the geometry of those sets, which is of importance
for deriving complexity results, has yet to be investigated. Moreover, the
linear equality constrained case~\cite{SGratton_CWRoyer_LNVicente_ZZhang_2019}
suggests that polling directions can be generated directly within the
manifold and the tangent spaces, thereby avoiding the need for projections.

In this paper, we study a class of direct-search methods dedicated to
optimization on (smooth) Riemannian manifolds. We first provide complexity
guarantees for these methods, which \rrev{are} complementary to the existing
asymptotic convergence theory for Riemannian direct search. We then address the
generation of polling sets that meet the requirements of our complexity analysis.
To this end, we extend the Euclidean notion of the cosine measure to embedded
submanifolds and propose two constructions of positive spanning sets for a
manifold. Our first construction projects Euclidean PSSs onto tangent spaces,
and is similar to that used by \rrev{Kungurtsev} et
al.~\cite{VKungurtsev_FRinaldi_DZeffiro_2024}. Our second construction builds
PSSs intrinsically in tangent spaces, and bears similarity with known approaches
for linear equality constraints~\cite{SGratton_CWRoyer_LNVicente_ZZhang_2019,
	LLiu_XZhang_2006}. Both strategies inherit geometric guarantees from their
Euclidean counterparts, which we quantify using both the cosine measure and
a new quantity called complexity measure, \rrev{which} arises from complexity results
of direct search. In the specific case of the sphere manifold, we show that
intrinsic PSSs exhibit provably better guarantees than projected PSSs. Our
numerical experiments validate our theoretical findings on Riemannian
optimization problems, and further illustrate the benefit of intrinsic
approaches. The code for reproducing our experiments can be found on
Github\footnote{
	\url{https://github.com/bastiencavarretta/intrinsic-vs-projected-directsearch}.}.

The rest of the paper is organized as follows.
Section~\ref{sec:dsalgo} recalls key notions from Riemannian optimization and
presents our direct-search framework. We then derive complexity results for this
method in Section~\ref{sec:complexity} assuming that polling directions are of
sufficient quality. Section~\ref{sec:dsdirs} formalizes the notion of good
directions by introducing a Riemannian variant of the cosine measure, and also
describes two ways of producing polling sets in tangent spaces.
In Section~\ref{sec:study_cm_sphere}, we estimate the cosine measure of both
strategies in the case of the sphere manifold, and argue for the intrinsic
design of polling directions.
Section~\ref{sec:num} compares direct-search schemes based on intrinsic
and projected PSSs on manifold optimization problems.
Section~\ref{sec:conc} summarizes our findings and offers several perspectives
on our work.

\paragraph{Notations:} $\llbracket a,b \rrbracket = \{a,\dots,b\}$ is the set
of integers between $a$ and $b$. Throughout our paper, $m$ and $n$ are positive
integers that correspond to the dimension of a manifold $\M$ and the dimension
of an ambient space $E$, respectively. In the Euclidean space $\R^n$, we
denote by $\langle u , v \rangle$ the canonical inner product between vectors
$u$ and $v$ \rev{and $\|u\|_2$ the associated Euclidean norm. We let
	$\{ e_1, \dots, e_n \}$ denote the coordinate vectors in $\R^n$, while the
	identity matrix $I_n \in \R^{n \times n}$} is the matrix whose columns are the
vectors $e_i$ in order. Given a Euclidean vector space
$(E,\langle \cdot,\cdot \rangle_E)$, the orthogonal projection onto a subspace
$F \subset E$ with respect to $\langle \cdot, \cdot \rangle_E$ is written
$\proj_F$.

\section{Riemannian direct-search framework}
\label{sec:dsalgo}

\subsection{Background on Riemannian optimization}
\label{ssec:background}

This section recalls a few standard notions from differential geometry that
are used in Riemannian optimization. We refer the reader to standard
references~\cite{PAAbsil_RMahony_RSepulchre_2008,NBoumal_2023} for a
comprehensive treatment.

Let $\M$ be an $m$-dimensional smooth manifold equipped with a Riemannian
metric $\langle\cdot,\cdot\rangle_x$, which defines an inner product on each
tangent space $\T_x\M$. The associated norm is denoted by
$\|v\|_x = \sqrt{\langle v,v\rangle_x}$ (often written $\norm{v}$ when the
reference point is clear from context).
Any family $\BM = (\B_x)_{x \in \M}$ such that $\B_x$ is an \rrev{orthonormal} basis
of $\T_x\M$ (in the sense of $\|v\|_x$) for any $x \in \M$ is called an
\emph{\rrev{orthonormal} basis adapted to the manifold $\M$}.
The Riemannian metric induces a notion of gradient for smooth functions
$f\colon\M\to\R$: the \emph{Riemannian gradient} $\grad f(x)\in \T_x\M$ is
the unique tangent vector satisfying
\begin{align*}
	\mathrm{D} f(x)[v] & =	\inner{v}{\grad f(x)}_x\,
	\text{ for all } v\in \T_x\M,
\end{align*}
where $\mathrm{D} f(x)[v]$ denotes the directional derivative of $f$
along $v$.

To move on the manifold along a tangent direction, Riemannian optimization
algorithms use a \emph{retraction}, i.e. a smooth mapping
\begin{equation*}
	\Retr_x\colon \T_x\M \to \M,
\end{equation*}
such that each curve $c(t)=\Retr_x(tv)$ satisfies $c(0)=x$ and $c'(0)=v$.
Intuitively, $\Retr_x(v)$ generalizes the Euclidean step $x+v$ while
ensuring that the iterate remains on the manifold. For many manifolds of
interest, simple global retractions
exist~\cite[Chapter 4]{PAAbsil_RMahony_RSepulchre_2008}; otherwise, they
are locally defined in a ball of radius $\varrho(x)>0$ around the origin
of $\T_x\M$. In that case, the length of the step at $x\in \M$ must be
limited to $\varrho(x)$.

It is common for the manifold $\M$ to be a subset of a Euclidean vector
space $\left(E,\langle \cdot, \cdot \rangle_E\right)$. In that case, and
provided $\M$ is endowed with the metric $\langle \cdot, \cdot \rangle_E$,
we say that $\M$ is an \emph{embedded submanifold} of $E$.
For any smooth function $\tilde f: E \to \R$, let $f$ be the restriction
of $\tilde f$ to $\M$. It follows that
\begin{equation} \label{eq:riem_eucl_gradient}
	\grad f(x) = \proj_x(\nabla \tilde f(x)),
\end{equation}
where $\nabla \tilde f(x)$ is the Euclidean gradient and
$\proj_x=\proj_{\T_x\M}$ denotes the orthogonal projection onto the
tangent space $\T_x\M$.
As an illustrative example, consider the unit sphere
$\mathbb{S}^{n-1} := \{x\in\R^n : \rev{\|x\|_2^2 = 1}\}$.
Here $\T_x\mathbb{S}^{n-1} = \{v\in\R^n : \langle v,x\rangle = 0\}$,
the projection is
\rev{$\proj_x=\proj_{T_x\mathbb{S}^{n-1}}: y \mapsto (I_n - xx^\top)y$},
and a common retraction is $\Retr_x(v) = \frac{x+v}{\rev{\|x+v\|_2}}$.

Finally, we will make use of parallel transport to combine or compare gradients
at different points of $\M$. Parallel transport moves a vector $v\in \T_x\M$
to another tangent space $\T_y\M$ along a curve connecting $x$ and $y$ while
preserving inner products~\cite[Section 10.3]{NBoumal_2023}.

\subsection{Riemannian direct-search algorithm}
\label{ssec:dsalgo}

In this section, we describe our Riemannian direct-search algorithm, that
bears close similarity with that studied in Kungurtsev et
al.~\cite[Algorithm 1]{VKungurtsev_FRinaldi_DZeffiro_2024}. We however point
out that our framework, given in Algorithm~\ref{algo:RDS}, focuses on
directional direct-search schemes~\cite{TGKolda_RMLewis_VTorczon_2003}. In
addition to being amenable to a complexity analysis, this class of algorithms
allows us to focus on the role played by the directions used throughout the
method.

\rev{At each iteration $k$, Algorithm~\ref{algo:RDS} generates a finite
	\emph{polling set} $D_k\subset \T_{x_k}\M$ of tangent directions, that
	represents a proxy for the (negative) Riemannian gradient.
	Each direction $d$ is associated with a retraction step
	$\Retr_{x_k}(\alpha_k d)$, where $\alpha_k>0$ is the current step size, and
	the algorithm then polls the objective at those retraction steps to find
	sufficient decrease~\eqref{eq:suffdec} in the objective value. If one
	retraction step satisfies this requirement (successful iteration), the iterate
	is updated and the step size may increase to better explore the variable space.
	Otherwise, the iterate does not change and the step size is reduced to seek
	improvement in a closer neighborhood of the current point.
	Note that the polling process (Step \ref{step3}) in
	Algorithm~\ref{algo:RDS} can be either \textit{opportunistic} or
	\textit{complete}. Opportunistic polling evaluates all directions in a
	sequential fashion, and stops as soon as a direction is found
	satisfying~\eqref{eq:suffdec}, moving to the next iteration without
	polling the remaining directions in $\D_k$. By contrast, complete polling
	evaluates all directions in $\D_k$ and selects the one achieving the
	largest decrease (if any). Both strategies yield the same evaluation
	cost on unsuccessful iterations.}

\begin{algorithm}
	\DontPrintSemicolon 
	\KwIn{$x_0\in \calM$, $\alphamax>0$, $\alpha_0\in (0,\alphamax], c>0,0<\gammadec<1<\gammainc$\rev{.}}
	\For{$k=0,1,\dots$} {
		Generate a \rev{finite} polling set $\D_k\subset \T_{x_k}\calM$ of tangent vectors\rev{.} \label{step2}\;
		If there exists $d\in \D_k$ such that
		\begin{equation}
			\label{eq:suffdec}
			f(\Retr_{x_k}(\alpha_k d))
			<
			f(x_k) - \frac{c}{2}\alpha_k^2 \norm{d}_{\rev{x_k}}^2,
		\end{equation}
		then declare iteration $k$ as \emph{successful}, and \rev{define}
		$$x_{k+1}:= \Retr_{x_k}( \alpha_k d)
			\text{ and }
			\alpha_{k+1}:=\min\{\gammainc\alpha_k,\alphamax\}.$$
		\label{step3} \;
		Otherwise, declare iteration $k$ as \emph{unsuccessful}, and \rev{define}
		$$x_{k+1}:= x_k \text{ and }\alpha_{k+1}:= \gammadec \alpha_k.$$
	}
	\caption{Riemannian direct-search framework}
	\label{algo:RDS}
\end{algorithm}

\rev{The theory of direct-search frameworks such as Algorithm~\ref{algo:RDS}
	relies on the sufficient decrease condition and the adaptive stepsize choice,
	that are both tied to backtracking line-search
	procedures~\cite[Section 3.7]{TGKolda_RMLewis_VTorczon_2003}. The choice of
	polling directions is also crucial, particularly for deriving complexity
	guarantees~\cite{KJDzahini_FRinaldi_CWRoyer_DZeffiro_2025}. Obtaining such
	results in the Riemannian setting is one of the goals of this paper, and thus
	we} identify desirable properties for polling sets by studying the complexity
of Algorithm~\ref{algo:RDS} in the next section.

\section{Complexity analysis}
\label{sec:complexity}

In this section, we derive a worst-case complexity bound for
Algorithm~\ref{algo:RDS}. Our goal is to bound the number of function
evaluations required to obtain a point whose Riemannian gradient norm is below
a prescribed tolerance $\varepsilon>0$. We will show that this number is of
order $\mathcal{O}(\varepsilon^{-2})$, thus matching the Euclidean bound
derived by Vicente~\cite{LNVicente_2013} in terms of $\varepsilon$.
Our analysis follows the Euclidean setting. We first state our assumptions
as well as auxiliary lemmas in Section~\ref{ssec:assum}. Our main results are
proven in Section~\ref{ssec:thms}.

\subsection{Assumptions and intermediate lemma}
\label{ssec:assum}

\rev{In our analysis, we rely on standard assumptions on the objective function. Our first \rrev{assumption} ensures that the problem is not unbounded below.}

\begin{assumption}\label{assu:flow}
	There exists $\flow\in \R$ such that $f(x) \geq \flow$ for all $x\in \calM$.
\end{assumption}

\rev{We also require a Lipschitz condition on the Riemannian gradient of $f$,
	similarly to existing manifold optimization
	\rrev{analyses}~\cite{KJDzahini_FRinaldi_CWRoyer_DZeffiro_2025,NBoumal_2023}.}

\begin{assumption}\label{assu:lipschitz_gradient}
	The function $f\colon \M\to \R$ is continuously differentiable and
	\rev{the gradient of} the pullback $f\circ \Retr_x$ is $L$-Lipschitz
	continuous at every $x$ \prev{in} 
	\begin{equation*}
		\calL_f(x_0):= \left\{x\in \M\colon f(x) \leq f(x_0)\right\}.
	\end{equation*}
\end{assumption}

Note that Assumption~\ref{assu:lipschitz_gradient} holds in particular if the
retraction is the exponential map and the Riemannian gradient $\grad f$ is a
$L$-Lipschitz continuous vector field~\cite[Definition 10.44]{NBoumal_2023}.
Under Assumption~\ref{assu:lipschitz_gradient}, we have
\begin{equation}\label{eq:Lipgrad}
	f(\Retr_x(v)) \leq f(x) + \inner{\grad f(x)}{v}_x + \dfrac{L}{2}\norm{v}_x^2
\end{equation}
for any $x\in \calL_f(x_0)$ and all $v\in \T_x\M$ such that
$\norm{v}_x \leq \varrho(x)$, where we recall that $\varrho(x)$ denotes the
radius of a ball over which the retraction $\Retr_x$ is well defined.
\rev{Moreover, the property~\eqref{eq:Lipgrad} is satisfied whenever $\M$ is
	a compact embedded submanifold of a Euclidean space $E$, and $f$ admits a
	$\tilde L$-smooth extension
	$\tilde f: E \to \R$ \cite[Lemma 4]{NBoumal_PAAbsil_CCartis_2019}.}

We now make our main assumption of interest in this paper, that relates to
the polling directions used by Algorithm~\ref{algo:RDS}.

\begin{assumption}\label{assu:pss}
	There exist $\kappa>0$ and
	$0<\Dmin\leq\Dmax < \frac{\varrho_0}{\alphamax}$ such that,
	for every iteration $k$, the \rev{finite}
	polling set $\D_k\subset \T_{x_k}\M$ satisfies
	$\Dmin \le \norm{d}_{x_k} \le \Dmax$ for all $d \in D_k$ and
	\begin{align}
		\label{eq:assump_cm_gradient_inequality}
		\dfrac{\inner{-\grad f(x_k)}{d}_{x_k}}{\norm{\grad f(x_k)}_{x_k}\norm{d}_{x_k}}
		\geq \kappa
		\quad \text{ for some }d\in D_k,
	\end{align}
	when $\norm{\grad f(x_k)}_{x_k} \neq 0$.
\end{assumption}

Assumption~\ref{assu:pss} guarantees that at least one polling direction forms
an acute angle with the negative gradient, uniformly over $k$. As we will see
in Section~\ref{sec:dsdirs}, generating $\D_k$ as a positive spanning set of
$\T_{x_k}\M$ is a way of satisfying~\eqref{eq:assump_cm_gradient_inequality}.

\begin{remark}
	\label{rem:injradius}
	To ensure that every retraction step $\Retr_{x_k}(\alpha_k d)$ is well-defined,
	we assume that $\alphamax \Dmax$ does not exceed $\inf_{x\in \M} \varrho(x)$,
	which is positive if the injectivity radius $\varrho_0$ of the manifold is
	positive~\cite[Remark 2.2]{NBoumal_PAAbsil_CCartis_2019}. This assumption is
	common in the analysis of Riemannian optimization
	methods~\cite{AAgarwal_NBoumal_BBullins_CCartis_2021,
		NBoumal_PAAbsil_CCartis_2019}.
\end{remark}

Having stated our assumptions, we now provide two straightforward extensions
of standard Euclidean arguments to the Riemannian setting. Our first result
relates the step size at a given unsuccessful iteration with the Riemannian
gradient at the corresponding iterate. Although the proof tracks that of the
Euclidean case, we provide it after the lemma to highlight the role of
\rev{Assumption~\ref{assu:pss}}.

\begin{lemma} \label{lemma:failure_small_gradient}
	Under Assumptions~\ref{assu:lipschitz_gradient} and~\ref{assu:pss}, if
	iteration $k$ is unsuccessful, \rrev{i.e.,}
	\begin{equation}\label{eq:unsucessfulness}
		f(x_k) - \frac{c}{2}\alpha_k^2 \|d\|_{x_k}^2
		\leq
		f\left(\Retr_{x_k}(\alpha_k d)\right) \text{ for all } d\in \D_k,
	\end{equation}
	then
	\begin{align}\label{eq:bound_gradient_unsuccessful}
		\norm{\grad f(x_k)}_{x_k} \leq \dfrac{(L + c)\Dmax}{2\kappa}\alpha_k.
	\end{align}
\end{lemma}
\begin{proof}
	Throughout the proof, we omit the index of the norms and inner
	products for sake of readability. Since~\eqref{eq:bound_gradient_unsuccessful}
	is straightforward when $\norm{\grad f(x_k)}=0$, we assume for the
	rest of the proof that the gradient is nonzero.

	For any $d \in \D_k$,  we combine~\eqref{eq:unsucessfulness} with the
	Lipschitz property \eqref{eq:Lipgrad} to get
	\begin{align*}
		f(x_k) -  \frac{c}{2}\alpha_k^2 \norm{d}^2
		\leq
		f\left(\Retr_{x_k}(\alpha_k d)\right) \leq f(x_k)
		+ \inner{\grad f(x_k)}{\alpha_k d} + \frac{L}{2}\alpha_k^2\norm{d}^2.
	\end{align*}
	This gives
	\begin{align}
		\label{ineq:bounded_innerproduct}
		\dfrac{\inner{- \grad f(x_k)}{d}}{\norm{d}}\leq \frac{(L+c)}{2}\alpha_k\|d\|.
	\end{align}
	By Assumption~\ref{assu:pss}, there exists $d\in D_k$ such that
	\eqref{eq:assump_cm_gradient_inequality} holds.
	Combining~\eqref{eq:assump_cm_gradient_inequality}
	with~\eqref{ineq:bounded_innerproduct} gives
	\begin{align*}
		\kappa \|\grad f(x_k) \|\leq \dfrac{\inner{-\grad f(x_k)}{d}}{\|d\|}
		\leq
		\frac{(L+c)\Dmax}{2}\alpha_k,
	\end{align*}
	which is~\eqref{eq:bound_gradient_unsuccessful}.
\end{proof}

\subsection{Main results}
\label{ssec:thms}

In this section, we derive complexity bounds on our algorithm. We first
\rrev{show that not all iterations can be successful, and provide a bound on
	the first unsuccessful iteration index.}

\begin{proposition}\label{prop:k0}
	Under Assumptions~\ref{assu:flow}, \ref{assu:lipschitz_gradient}
	and \ref{assu:pss}, \rrev{there exists a least one unsuccessful iteration.
		Moreover, letting $k_0$ denote the index of the first unsuccessful
		iteration, we have}
	\begin{align}
		\label{eq:k0}
		k_0
		\leq \dfrac{2(f(x_0) - \flow)}{c \alpha_0^2 \Dmin^2}
		\rev{< \infty}.
	\end{align}
\end{proposition}
\begin{proof}
	We again omit indices in norms throughout the proof.

	\rrev{Suppose first that all iterations are unsuccessful. Then, for
		any $K \ge 1$, every iteration of index $k \in \llbracket 0,K-1\rrbracket$
		is successful.} Let $d_k \in \D_k$ denote the direction taken at iteration
	\rrev{$k\in \llbracket 0,K-1\rrbracket$}. By the sufficient decrease
	condition~\eqref{eq:suffdec}, we have
	\begin{align*}
		f(x_k) - f(x_{k+1})
		\geq \frac{c}{2} \alpha_k^2 \|d_k\|^2.
	\end{align*}
	Because the step size has not been reduced before the first unsuccessful
	iteration, $\alpha_k \geq \alpha_0$ for all $k < \rrev{K}$.
	Using $\|d_k\| \ge \Dmin$ yields
	\begin{equation}
		\label{eq:decsucck0}
		f(x_k) - f(x_{k+1})
		\ge \frac{c}{2}\,\alpha_0^2\,\Dmin^2.
	\end{equation}
	Summing~\eqref{eq:decsucck0} over \rrev{$k \in \llbracket 0, K-1 \rrbracket$}
	and applying Assumption~\ref{assu:flow} yields
	\begin{equation}
		\label{eq:boundk0_telescoping_argument}
		f(x_0) - \flow
		\ge \sum_{k = 0}^{\rrev{K} - 1} \big(f(x_k) - f(x_{k+1})\big)
		\ge \sum_{k = 0}^{\rrev{K}- 1} \,\frac{c}{2}\,\alpha_0^2\,\Dmin^2
		=\rrev{K}\,\frac{c}{2}\,\alpha_0^2\,\Dmin^2.
	\end{equation}
	\rrev{
		Since $f(x_0)-\flow$ is finite, we conclude that $K < \infty$. Thus
		not all iterations can be successful, and the index $k_0$ is well defined.
		Moreover, all iterations with index $k \in \llbracket 0,k_0-1 \rrbracket$
		are successful, and thus~\eqref{eq:boundk0_telescoping_argument} holds
		with $K=k_0$, yielding~\eqref{eq:k0}.}
\end{proof}

We then bound the number of successful iterations that follow the first
unsuccessful iteration.
\begin{theorem} \label{thm:bound_success}
	Under Assumptions~\ref{assu:flow}, \ref{assu:lipschitz_gradient}
	and \ref{assu:pss}, let $k_0$ be the index of the first unsuccessful
	iteration of Algorithm~\ref{algo:RDS}. Given any $\varepsilon \in (0,1)$,
	assume that $\|\grad f(x_{k_0})\|_{x_{k_0}}> \varepsilon$ and let $j_1$ be the first
	iteration after $k_0$ such that $\norm{\grad f(x_{j_1 + 1})}_{x_{j_1+1}}\leq \varepsilon$.
	Finally, let $|\Skj|$ be the number of successful iterations between
	$k_0$ and $j_1$ \rrev{($j_1$ excluded)}. Then,
	\begin{align}
		\label{eq:bound_success}
		|\Skj|
		\leq
		\dfrac{\left( f(x_0) - \flow\right)}{2c}\left(
		\dfrac{(L +c) \Dmax}{ \gammadec\Dmin }
		\right)^2 \kappa^{-2}\varepsilon^{-2}.
	\end{align}
\end{theorem}
\begin{proof}
	\rrev{Note that $j_1 > k_0$ by definition, since the former must be
		successful while the later is unsuccessful.}
	We consider indices \rrev{$k=k_0,\dots,j_1-1$, for which
		$\norm{\grad f(x_k)} > \varepsilon$}. If iteration $k$ is unsuccessful,
	Lemma~\ref{lemma:failure_small_gradient} gives
	\begin{align}
		\label{eq:unsucck0j1}
		\alpha_k
		\geq
		\dfrac{2\kappa}{(L + c)\Dmax}\|\grad f(x_k)\| > \dfrac{2\kappa}{(L + c)\Dmax} \varepsilon.
	\end{align}
	On the other hand, if $k$ corresponds to a successful iteration, then
	\begin{align}
		\label{eq:succk0j1}
		\alpha_k \geq \gammadec \alpha_{k_1},
	\end{align}
	where $k_1$ is the most recent unsuccessful iteration preceding $k$ (with
	possibly $k_1=k_0$). Combining~\eqref{eq:succk0j1} and~\eqref{eq:unsucck0j1}, we
	obtain
	\begin{align}\label{eq:lower_bound_stepsize}
		\alpha_k > \gammadec \dfrac{2\kappa}{(L + c)\Dmax} \varepsilon,
		\quad \text{for all } k = k_0, \dots, \rrev{j_1-1}.
	\end{align}

	Therefore, if iteration $k>k_0$ is successful and $x_{k+1}=\rev{\Retr_{x_k}(\alpha_k d)}$, the
	sufficient decrease condition~\eqref{eq:suffdec} implies
	\begin{equation}
		\label{eq:suffdeck0j1}
		f(x_k) - f(x_{k+1})
		\geq \frac{c}{2}\alpha_k^2 \norm{d}^2
		\geq \frac{c}{2} \gammadec^2 \left(\dfrac{2\kappa}{(L + c)\Dmax} \right)^2\varepsilon^2 \Dmin^2.
	\end{equation}
	Using Assumption~\ref{assu:flow} and summing~\eqref{eq:suffdeck0j1}
	\rrev{over successful iterations in $\llbracket k_0,j_1-1 \rrbracket$}
	leads to
	\begin{eqnarray*}
		f(x_0) - \flow
		\geq f(x_{k_0}) - f(x_{j_1})
		=
		\sum_{k =k_0}^{j_1-1} \rev{(f(x_k) - f(x_{k+1}))}
		&\rrev{=}
		&\rrev{\sum_{k \in \Skj} (f(x_k) - f(x_{k+1}))} \\
		&\ge
		&|\Skj|\frac{c}{2}  \left(
		\dfrac{2\kappa \gammadec \Dmin}{(L + c)\Dmax} \right)^2\varepsilon^2.
	\end{eqnarray*}
	Rearranging the above inequality, we arrive at~\eqref{eq:bound_success}.
\end{proof}

Our next step towards a complexity result consists in relating the number of
unsuccessful iterations to that of successful iterations. The proof of this
result is a consequence of the stepsize updating rule, and thus follows directly
from the Euclidean case.

\begin{theorem} \label{thm:bound_failure}
	Under the assumptions of Theorem~\ref{thm:bound_success},
	let $|\Ukj|$ denote the number of unsuccessful iterations of
	Algorithm~\ref{algo:RDS} between indices $k_0$ and $j_1$.

	Then,
	\begin{align}
		\label{eq:bound_failure}
		|\Ukj|
		\leq
		\left|\dfrac{\log(\gammainc)}{\log(\gammadec)}\right| |\Skj| - \dfrac{\log(\alpha_{k_0})}{\log(\gammadec)}  +\dfrac{\log\left( \gammadec \dfrac{2\kappa}{(L + c)\Dmax} \varepsilon \right)}{\log(\gammadec)}.
	\end{align}
\end{theorem}
\begin{proof}
	At each iteration, the step size update gives
	$\alpha_{k+1} = \gammadec\,\alpha_k$ if iteration $k$ is unsuccessful,
	or $\alpha_{k+1} \leq \gammainc\alpha_k$ if iteration $k$ is successful.
	Applying this rule inductively from $k_0$ to $j_1$ gives
	\begin{align}
		\label{eq:intermediatebound_failure}
		\alpha_{j_1} \leq \alpha_{k_0} \gammadec^{|\Ukj|}\gammainc^{|\Skj|} \rev{.}
	\end{align}
	Taking logarithms on both sides of~\eqref{eq:intermediatebound_failure}
	and dividing by $\log(\gammadec)<0$ yields
	\begin{align}\label{eq:bound_Ukj_proof}
		|\Ukj|
		\leq
		\left|\dfrac{\log(\gammainc)}{\log(\gammadec)}\right| |\Skj|
		- \dfrac{\log(\alpha_{k_0})}{\log(\gammadec)}
		+ \dfrac{\log(\alpha_{j_1})}{\log(\gammadec)}.
	\end{align}
	To bound the last term on the right-hand side, we \rrev{note that the lower bound
		from~\eqref{eq:lower_bound_stepsize} holds for $k = j_1$, implying that}
	\begin{align}
		\label{eq:lower_bound_stepsize_Ukj}
		\alpha_{j_1}
		> \gammadec \frac{2\kappa}{(L + c)\Dmax}\,\varepsilon.
	\end{align}
	Plugging~\eqref{eq:lower_bound_stepsize_Ukj} into~\eqref{eq:bound_Ukj_proof}
	completes the proof.
\end{proof}

We can now establish worst-case complexity bounds for Algorithm~\ref{algo:RDS}
by combining the bounds above.

\begin{theorem} \label{theorem:evaluation_complexity}
	\rev{Let Assumptions~\ref{assu:flow}, \ref{assu:lipschitz_gradient} and
		\ref{assu:pss} hold. For any $\varepsilon \in (0,1)$,}
	Algorithm~\ref{algo:RDS} reaches a point $x \in \M $ satisfying
	$\|\grad f(x)\|_x\le \varepsilon$ in at most
	\begin{align} \label{eq:iteration_upperbound}
		C\kappa^{-2}\varepsilon^{-2} & ~~\text{iterations},
	\end{align}
	and
	\begin{align}\label{eq:evaluation_upperbound}
		C\cardpss \kappa^{-2}\varepsilon^{-2} & ~~\text{function evaluations,}
	\end{align}
	where $\kappa$ is defined in Assumption \ref{assu:pss},
	$r$ is an upper bound on the cardinality of the polling sets, and $C$ is a
	positive constant that depends on
	$\gammadec,~\gammainc,~\alpha_0,~\alphamax,~\Dmin,~\Dmax,~c,~L,~f(x_0),~\flow$.
\end{theorem}
\begin{proof}
	\rev{
	Let $k_0$ be the index of the first unsuccessful iteration, which
	exists by Proposition~\ref{prop:k0}.

	Suppose first that $\|\grad f(x_{k_0})\|_{x_{k_0}} > \varepsilon$, and let
	$j_1$ be defined as in Theorem~\ref{thm:bound_success}. Combining the bounds
	obtained for $k_0$ (Proposition~\ref{prop:k0}), $|\Skj|$
	(Theorem~\ref{thm:bound_success}) and $|\Ukj|$
	(Theorem~\ref{thm:bound_failure}) then gives
	\begin{eqnarray}
		\label{eq:proof_evaluation_complexity1}
		j_1 \rrev{+1}
		&\rrev{=} &k_0 + |\Skj| + |\Ukj| \rrev{+1} \nonumber \\
		&\le &\dfrac{2(f(x_0) - \flow)}{c \alpha_0^2 \Dmin^2}
		+\left (1+ \rrev{\left | \frac{\log \gammainc}{\log \gammadec}\right |}\right )
		\frac{f(x_0)-\flow}{2c}\left ( \frac{(L+c)\Dmax}{\gammadec \Dmin}\right )^2
		\kappa^{-2} \varepsilon^{-2} \nonumber \\
		& &- \frac{\log\alpha_{k_0}}{\log \gammadec}
		+\frac{\log\left(
			\gammadec \frac{2\kappa}{(L + c)\Dmax} \varepsilon
			\right)}{\log(\gammadec)} \rrev{+1}.
	\end{eqnarray}
	\rrev{Using that $\alpha_{k_0}\le \alphamax$ and defining}
	$
		C_1:=\left (
		1 + \rrev{\left | \frac{\log \gammainc}{\log \gammadec}\right |} \right ) \frac{f(x_0)-\flow}{2c}
		\left (\frac{(L+c)\Dmax}{\gammadec \Dmin}\right )^2,
	$
	as well as
	$
		C_2:=- \frac{\log\alphamax}{\log \gammadec}+
		\frac{\log\left(  \frac{2\gammadec}{(L + c)\Dmax} \right)}{\log(\gammadec)}
		+\dfrac{2(f(x_0) - \flow)}{c \alpha_0^2 \Dmin^2}\rrev{+1},
	$
	the bound~\eqref{eq:proof_evaluation_complexity1} becomes
	\begin{equation*}
		j_1 \rrev{+1 }
		\le
		C_1 \kappa^{-2} \varepsilon^{-2}
		+ \frac{\log(\kappa \varepsilon)}{\log \gammadec}
		+ C_2
		\leq
		C_1 \kappa^{-2} \varepsilon^{-2}
		+ \frac{|\log(\kappa \varepsilon)|}{|\log \gammadec|}
		+ C_3\times 1
	\end{equation*}
	with $C_3 :=
		\rrev{\left | \frac{\log\alphamax}{\log \gammadec} \right |}
		+ \left |\frac{\log\left(  \frac{2\gammadec}{(L + c)\Dmax} \right)}
		{\log(\gammadec)}\right |
		+ \left | \dfrac{2(f(x_0) - \flow)}{c \alpha_0^2 \Dmin^2}\right | \rrev{+1}
		\ge |C_2|$.
	Using that $|\log(x)| \leq x^{-2}$ and $x^{-2}>1$ for all $x \in (0,1)$
	yields
	\rrev{
		\begin{equation*}
			j_1  \rrev{+1}
			\leq \left ( C_1+\frac{1}{|\log \gammadec|}
			+ C_3\right )
			\kappa^{-2} \varepsilon^{-2}
			= C \kappa^{-2} \varepsilon^{-2},
		\end{equation*}}
	where $C := \left ( C_1+\frac{1}{|\log \gammadec|}+ C_3\right )$.
	Thus, \eqref{eq:iteration_upperbound} holds.
	}

	\rev{Suppose now that $\|\grad f(x_{k_0})\|_{x_{k_0}} \le \varepsilon$.
	Then, the number of iterations to reach an $\varepsilon$-stationary
	point is at most $k_0$. Applying Proposition~\ref{prop:k0} gives
	\[
		k_0 \leq \dfrac{2(f(x_0) - \flow)}{c \alpha_0^2 \Dmin^2}
		\leq C_3 \leq C_3 \kappa^{-2} \varepsilon^{-2}
		\leq C \kappa^{-2} \varepsilon^{-2},
	\]
	proving that~\eqref{eq:iteration_upperbound} also holds in that case.}

	\rev{Finally, the bound~\eqref{eq:evaluation_upperbound} follows
		directly from~\eqref{eq:iteration_upperbound} since each iteration
		requires at most $r$ function evaluations.}
\end{proof}

The complexity bound established in Theorem~\ref{theorem:evaluation_complexity} matches
that of Euclidean direct-search methods~\cite{KJDzahini_FRinaldi_CWRoyer_DZeffiro_2025}.
In particular, its dependencies on $r$ and on the parameter $\kappa$, that capture the
size and the geometry of the polling directions, are identical.  In the Euclidean case
$\M=\R^m$, directions can be chosen so that \rrev{$r\kappa^{-2}$} exhibits a polynomial
\rrev{dependence} in $m$, with the best \rrev{possible dependence} being
$\mathcal{O}(m^2)$~\cite{MDodangeh_LNVicente_ZZhang_2016}.
In the next section, we will explore constructions of polling sets that allow to
quantify both $r$ and $\kappa$ in a Riemannian setting.

\section{Choosing poll directions on a Riemannian manifold}
\label{sec:dsdirs}

In this section, we construct appropriate polling sets for Riemannian
optimization. \rev{In the Euclidean case, polling directions are typically
	chosen to form finite positive spanning sets (PSSs), uniformly controlling
	their geometric quality as well as their cardinality.}
We recall the key definitions and examples associated with PSSs in
Section~\ref{ssec:PSSs_Rn}, then generalize those \rev{concepts} to the
Riemannian setting in Section~\ref{ssec:PSSs_TxM}. Afterwards, we present two
approaches for generating PSSs in a Riemannian sense from a PSS in \rev{a}
Euclidean sense.

\subsection{Positive spanning sets in \texorpdfstring{$\R^m$}{Rm}}
\label{ssec:PSSs_Rn}

The concept of positive spanning set is closely related to the notion of
cosine measure, that has played a major role in the analysis of direct-search
methods~\cite{TGKolda_RMLewis_VTorczon_2003}. Since both notions are of
interest for this paper, we define PSSs using the cosine measure as follows.

\begin{definition} \label{def:cm_Rn}
	Let $\D$ be a \rev{finite} set of vectors in $\R^m\setminus \{0\}$
	\rev{with $m\geq 1$}.
	The \emph{cosine measure of $\D$} is defined as
	\begin{equation}
		\cm(\D):= \underset{\substack{v\in \R^m \\ v \neq 0}}{\min} \,
		\underset{d\in \D}{\max}\dfrac{\inner{d}{v}}{\|d\|\|v\|}.
	\end{equation}
	A set $\D$ such that $\cm(\D)>0$ is called a \emph{positive spanning set}
	of $\R^m$.
\end{definition}

Note that Definition~\ref{def:cm_Rn} can be extended to general
Euclidean spaces. However, we focus on the case of $\R^m$ both for
simplicity and consistency with our experiments in subsequent sections.

We present below three examples of positive spanning sets that can be
obtained from any \rrev{orthonormal} basis of $\R^m$, that will be used later to
design polling sets in the Riemannian setting. Arguably the most classical
construction consists in using the vectors in the basis and their
negatives.

\begin{example}[Positive and negative basis\rev{\cite{GNaevdal_2019}}]
	\label{example:pssa}
	Let $B=[b_1,\dots,b_m] \in \R^{m \times m}$ be an orthogonal matrix.
	The set
	\begin{equation}
		\label{eq:pssa}
		\pssa(B) := \left\{b_1,\dots,b_m,{-}b_1,\dots,{-}b_m \right\}
	\end{equation}
	is a positive spanning set of $\R^m$ with $2m$ elements and
	cosine measure $\tfrac{1}{\sqrt{m}}$.
\end{example}

When $B=I_m$, the construction~\eqref{eq:pssa} corresponds to using
coordinate vectors and their negatives. Using those as polling directions
gives rise to the coordinate search method~\cite{CAudet_WHare_2017}.

A PSS with smaller cardinality can be constructed from an \rrev{orthonormal}
basis by using the sum of negative vectors.

\begin{example}[Basis and sum of negatives] \label{example:pssb}
	Let $B=[b_1,\dots,b_m] \in \R^{m \times m}$ be an orthogonal matrix.
	The set
	\begin{equation}
		\label{eq:pssb}
		\pssb(B) := \left\{b_1,\dots,b_m,
		{-}\frac{1}{\sqrt{m}}\sum_{i=1}^m b_i \right\}
	\end{equation}
	is a positive spanning set of $\R^m$ with $m+1$ elements and
	cosine measure $\tfrac{1}{\sqrt{m^2 + 2(m-1)\sqrt{m}}}$.
\end{example}

The cosine measure of $\pssb(B)$ was only recently
\rev{derived~\cite{WHare_GJarryBolduc_SKerleau_CWRoyer_2024,GJarryBolduc_2023}}. It can be shown that obtaining the largest
cosine measure using $m+1$ vectors requires a more involved
construction, which we describe below.

\begin{example}[Uniform angles~\cite{ARConn_KScheinberg_LNVicente_2009}]
	\label{example:pssc}
	Let $G_m \in \R^{m \times m}$ be the matrix with $1$s on its
	diagonal and $-\tfrac{1}{m}$ elsewhere. Let $G_m = L_m\,L_m^\T$
	be the Cholesky decomposition of $G_m$ (i.e. $L_m$ is lower
	triangular with positive diagonal elements), and denote by
	$v_1,\dots,v_m$ the columns of $L_m^\T$.

	For any orthogonal matrix $B \in \R^{m \times m}$,
	the set
	\begin{equation}
		\label{eq:pssc}
		\pssc(B) := \left\{ B v_1,\dots,B v_m,-\sum_{i=1}^m B v_i \right\}
	\end{equation}
	is a PSS with $m+1$ vectors and cosine measure $\frac{1}{m}$.
\end{example}

Note that given two orthogonal matrices $B=[b_1,\dots,b_m]$ and $\tilde{B}$ in
$\R^{m \times m}$, we have
\[
	\pssa(\tilde{B})
	=
	\left\{ \pm \tilde{B}\,B^\T b_i\ \middle| i=1,\dots,m \right\}.
\]
Similar formulas hold for the other PSS choices. As a result, we
will omit the basis $B$ when it is clear from context. In particular,
we will use $\pssa$, $\pssb$, $\pssc$ to denote $\pssa(I_m)$, $\pssb(I_m)$
and $\pssc(I_m)$, respectively.

\begin{table}[htb!]
	\centering
	\caption{Common PSSs in $\R^m$ and their properties}
	\label{table:usualEuclideanPSSs}
	\begin{tabular}{lccc}
		\toprule
		PSS                       & $\pssa$              & $\pssb$                                 & $\pssc$            \\
		\midrule
		Example                   & \ref{example:pssa}   & \ref{example:pssb}                      & \ref{example:pssc} \\

		\midrule
		representation in $\R^2$  &
		\begin{tikzpicture}[scale = 0.6, baseline={(current bounding box.center)}]
			\draw[->,thick] (0,0) -- (0,1);
			\draw[->,thick] (0,0) -- (0,-1);
			\draw[->,thick] (0,0) -- (1,0);
			\draw[->,thick] (0,0) -- (-1,0);
		\end{tikzpicture}
		                          &
		\begin{tikzpicture}[scale = 0.6, baseline={(current bounding box.center)}]
			\draw[->,thick] (0,0) -- (0,1);
			\draw[->,thick] (0,0) -- (1,0);
			\draw[->,thick] (0,0) -- (-0.707,-0.707);
		\end{tikzpicture}
		                          &
		\begin{tikzpicture}[scale = 0.6, baseline={(current bounding box.center)}]
			\begin{scope}[rotate = 0]
				\draw[->,thick] (0,0) -- (1,0);
				\draw[->,thick] (0,0) -- (-1/2,0.866);
				\draw[->,thick] (0,0) -- (-1/2,-0.866);
			\end{scope}
		\end{tikzpicture}                                                 \\
		\midrule
		cardinality               & $2m$                 & $m+1$                                   & $m+1$              \\
		\midrule
		cosine measure $\cm$      & $\frac{1}{\sqrt{m}}$ & $\frac{1}{\sqrt{m^2+ 2(m-1) \sqrt{m}}}$
		                          & $\frac{1}{m}$                                                                       \\
		\midrule
		complexity measure $\chi$ & $2m^2$               & $m^3 + O(\rev{m^{5/2}})$                & $m^3 + O(m^2)$     \\
		\bottomrule
	\end{tabular}
\end{table}

We end this section by introducing the \emph{complexity measure} of a PSS,

\begin{definition}\label{def:complexity_measure_Rn}
	Let $\D$ be a positive spanning set of $\R^m$. The value
	\begin{align} \label{eq:complexity_measure_Rn}
		\chi(\D):= |\D|\cm(\D)^{-2}
	\end{align}
	is called the \emph{complexity measure of $\D$}.
\end{definition}

Table~\ref{table:usualEuclideanPSSs} reports complexity measures for the three PSSs of
interest $\pssa$, $\pssb$ and $\pssc$. In particular, one observes \rev{that $\pssa$ corresponds
	to} the complexity measure with mildest dependency on $m$. This dependency is
known to be optimal in order~\cite{MDodangeh_LNVicente_ZZhang_2016}.

The complexity measure naturally arises in complexity bounds for direct
search~\cite[Theorem 2.3]{KJDzahini_FRinaldi_CWRoyer_DZeffiro_2025}. We
further emphasize the connection with our analysis from
Section~\ref{sec:complexity} below.

\begin{corollary}
	\label{coro:complexitychi}
	Let the assumptions of Theorem~\ref{theorem:evaluation_complexity} hold
	for some $\varepsilon \in (0,1)$. Suppose further that $\M=\R^m$, and
	that Algorithm~\ref{algo:RDS} is implemented using the same PSS $\D$ of $\R^m$
	at every iteration, i.e. $\D_k = \D$ for every $k$.
	Then, Algorithm~\ref{algo:RDS} reaches a point
	$x \in \R^m $ satisfying $\|\grad f(x)\|\le \varepsilon$ in at most
	\begin{align}
		\label{eq:complexitychi}
		C\,\chi(\D)\varepsilon^{-2}
	\end{align}
	function evaluations, where $C>0$ is the constant defined in
	Theorem~\ref{theorem:evaluation_complexity}.
\end{corollary}

\rev{
	We point out that bounds such as~\eqref{eq:complexitychi} are those actually
	used to assess how the number of function evaluations is affected by the
	choice of polling directions in the worst case. In particular,
	Corollary~\ref{coro:complexitychi} is established for a method that uses the
	same PSS at every iteration (akin to coordinate search), but the result also
	applies when using random rotations of the same PSS at every iteration. In
	addition, variants relying on opportunistic and complete polling share the
	same complexity properties, as the latter corresponds to a worst-case
	scenario of the former in terms of function evaluations used per iteration.}

\subsection{Positive spanning sets in \texorpdfstring{$\M$}{M}} \label{ssec:PSSs_TxM}

In this section, we aim at generalizing the notions of PSS, cosine measure and complexity
measure to the Riemannian setting, bearing in mind that those PSSs will be used for polling
within Algorithm~\ref{algo:RDS}.

In that context, a natural strategy consists in using positive spanning sets for each tangent
space encountered by the algorithm. This approach was used by Kungurtsev et al.~\cite{VKungurtsev_FRinaldi_DZeffiro_2024}, and bears connection with earlier work by
Dreisigmeyer~\cite{DWDreisigmeyer_2007a}. We formalize below the concept of PSS used in these
works.

\begin{definition}\label{def:cm_TxM}
	\rev{Let $x \in \M$ be a point in a manifold with positive dimension. Let $\D_x$} be a \rev{finite} set of tangent vectors in the tangent space $\T_x\M$.
	The \emph{cosine measure of $\D_x$} is defined as
	\begin{equation}
		\label{eq:cm_TxM}
		\cmx(\D_x):=
		\underset{\substack{v\in \T_x\M \\ v \neq 0}}{\min} \,
		\underset{d\in \D_x}{\max}\dfrac{\inner{d}{v}_x}{\|d\|_x\|v\|_x}.
	\end{equation}
	A set $\D_x$ such that $\cmx(\D_x)>0$ is called a positive spanning set (PSS) of $\T_x\M$.
\end{definition}

\rev{The cosine measure~\eqref{eq:cm_TxM} is a special case of the cosine
	measure with respect to a subspace~\cite{CAudet_WHare_GJarryBolduc_2025} where
	all vectors lie within the subspace. However, in a manifold setting, the reference
	point $x$ changes over the course of the algorithm, which changes the subspace as
	well. As a result, Definition~\ref{def:cm_TxM} alone is not sufficient for our
	purpose. We propose the construction below as a way to account for the entire
	manifold.}

\begin{definition}\label{def:cm_wholeM}
	Let $\DM = (\D_x)_{x \in \M}$ be a family of sets, where $\D_x$ is a \rrev{finite} set of vectors in
	$\T_x\M$ for every $x \in \M$.
	The \emph{cosine measure of family $\DM$} is defined as
	\begin{equation} \label{eq:cosine_measure_wholeM}
		\cm_\M(\DM) := \inf_{x \in \M} \cmx(\D_x).
	\end{equation}
	A family $\DM$ such that $\cm_{\M}(\DM)>0$ is called a positive spanning set (PSS) of
	$\M$.
\end{definition}

When $\M$ \prev{is the Euclidean space $\R^m$ equipped} with its canonical metric, any
tangent space $\T_x\M$ is equal to $\R^m$. The previous definition thus adequately
generalizes Definition \ref{def:cm_Rn}, while accounting for the fact that a
different PSS may be required for each tangent space in general. Provided one can
generate a PSS $\DM$ for $\M$ with $\cm_\M(\DM) \ge \kappa$, Assumption~\ref{assu:pss}
will hold and the method will be endowed with theoretical guarantees. To better
quantify those guarantees according to the PSS at hand, we propose the following
generalization of the complexity measure.

\begin{definition} \label{def:complexity_measure_M}
	Let $\DM := \left (\D_x \right ) _{x \in \M}$ be a PSS of $\M$. The value
	\begin{align}
		\label{eq:complexity_measure_M}
		\chi_\M(\DM):= \left ( \sup_{x \in \M} |\D_x|\right ) \cm_\M(\DM)^{-2}
	\end{align}
	is called the \emph{complexity measure of $\DM$}.
\end{definition}

When $\DM = (\D)_{x \in \R^m}$ is a family of identical PSSs of the manifold
$\R^m$, then $\chi_{\R^m}(\DM) = \chi(\D)$, recovering the complexity measure
defined for $\R^m$. A closer look at the analysis of Section~\ref{sec:complexity}
reveals that complexity bounds can be obtained using the value $\chi_\M(\DM)$
instead of the factor $r \kappa^{-2}$.

\begin{corollary}
	\label{coro:complexitychiM}
	Let the assumptions of Theorem~\ref{theorem:evaluation_complexity} hold
	for some $\varepsilon \in (0,1)$. Suppose further that Algorithm~\ref{algo:RDS}
	is implemented using a PSS $\DM = (\D_x)_{x \in \rev{\M}}$ of $\M$ so that
	$\D_k = \D_{x_k}$ for every $k$. Then, Algorithm~\ref{algo:RDS} reaches a point
	$x \in \M $ satisfying $\|\grad f(x)\|\le \varepsilon$ in at most
	\begin{align}
		\label{eq:complexitychiM}
		C\,\chi_{\M}(\DM)\varepsilon^{-2}
	\end{align}
	function evaluations, where $C>0$ is the constant defined in
	Theorem~\ref{theorem:evaluation_complexity}.
\end{corollary}

In the \rev{remainder} of Section~\ref{sec:dsdirs}, we will propose two ways
of generating positive spanning sets for $\M$ with provable guarantees
on their complexity measure.

\subsection{Building PSSs for \texorpdfstring{$\M$}{M} in an intrinsic fashion} \label{ssec:intrinsic_PSS}

Our first PSS construction aims at adapting the PSSs presented in Section~\ref{ssec:PSSs_Rn} to
tangent spaces. Given $x \in \M$, suppose that $\B_x=(b_{x,i})_{i=1,\dots,m}$ is an \rrev{orthonormal}
basis of $\T_x\M$. The basis $\B_x$ induces an isometric isomorphism
$\isomintr{\B_x}: \R^m \to \T_x\M$ such that $\isomintr{\B_x}(e_i)=b_{x,i}$ for every
$i=1,\dots,m$. Therefore, we can map any PSS $\D$ of $\R^m$ to the set
\begin{equation}
	\label{eq:intrinsicTxm}
	\D(\B_x) := \isomintr{\B_x}(\D) := \left\{\isomintr{\B_x}(d),~ d \in \D \right\}.
\end{equation}
Since \rrev{$\isomintr{\B_x}$} is an isometric isomorphism, the sets $\D$ and $\D(\B_x)$ have the same
cosine measure and cardinality. It follows that $\D(\B_x)$ is a PSS of $\T_x\M$, which we
call the \emph{intrinsic PSS induced by $\B_x$ and $\D$}. Applying this construction to the
entire manifold yields the following definition.

\begin{definition}[Intrinsic PSS] \label{def:intrinsic_PSS}
	Let $\BM = (\B_x)_{x \in \M}$ be a family of \rrev{orthonormal} bases such that $\B_x$ is an
	\rrev{orthonormal} basis for $\T_x\M$ and let $(\isomintr{\B_x})_{x \in \M}$ be the corresponding
	family of isometric isomorphisms. Let $\D$ be a PSS of $\R^m$. Then, the family
	\begin{equation}
		\label{eq:intrinsicM}
		\D(\BM) = (\D(\B_x))_{x \in \M},
	\end{equation}
	where $\D(\B_x)$ is given by~\eqref{eq:intrinsicTxm}\prev{,} defines a PSS of $\M$ called the
	\emph{intrinsic PSS} induced by $\BM$ and $\D$.
\end{definition}

With Definition~\ref{def:intrinsic_PSS} at hand, one can construct PSSs for $\M$ based on
the standard PSSs $\pssa$, $\pssb$, $\pssc$. For instance, given a family $\BM$, the
set $\pssaintr{\BM}$ has cosine measure $\tfrac{1}{\sqrt{m}}$ and complexity measure
$2m^2$, corresponding to the values for $\pssa$ as a PSS of $\R^m$.

\begin{remark}
	\label{rem:paralleltptisom}
	Since parallel transport is an isometric isomorphism between tangent spaces,
	one may consider defining PSS using parallel transport from a given
	$x_0 \in \M$ and $\D_{x_0} \in \T_{x_0}\M$. However, parallel transport is not
	a transitive mapping, and in particular is only well-defined within the injectivity
	radius~\cite[Chapters 4 and 7]{JMLee_2018}. Defining PSSs using parallel
	transport thus seems less straightforward \rev{than} our approach.
\end{remark}

\prev{The geometric properties of intrinsic PSSs are governed by the dimension of the manifold ($m$) rather than that of the ambient space ($n$), as they are defined independently of this ambient space.}
However, note that constructing an intrinsic PSS for $\M$ requires \rrev{computing} a family $\BM$
of \rrev{orthonormal} bases adapted to the manifold, which is not necessarily \prev{straightforward.} The next section investigates another construction that leverages knowledge about
the ambient space to alleviate the need for \rrev{orthonormal} bases.

\subsection{Building PSSs for \texorpdfstring{$\M$}{M} by projecting from the ambient space} \label{ssec:projected_PSS}

In this section, we suppose that the manifold $\M$ is an embedded submanifold \rev{of positive dimension in the Euclidean
	vector space $\R^n$}. Our approach will follow that proposed by Kungurtsev et
al.~\cite{VKungurtsev_FRinaldi_DZeffiro_2024} for such embedded submanifolds, and reduces to
known techniques when $\M$ is a linear subspace of
$\R^n$~\cite{CAudet_WHare_GJarryBolduc_2025,SGratton_CWRoyer_LNVicente_ZZhang_2019,
	LLiu_XZhang_2006}.

\rev{
	\begin{proposition}
		\label{prop:projPSSTxM}
		Let $\D$ be a PSS of $\R^n$. Let $x \in \M$, where $\M$ is an embedded submanifold
		in $\R^n$ with positive dimension. \prev{Let $\proj_x$ denotes the projection from
			$\R^n$ to $\T_x\M$.}  Then, the set
		\begin{equation}
			\label{eq:projPSSTxM}
			\fP_x(\D):=
			\left \{\frac{\proj_x(d)}{\|\proj_x(d)\|_x}~ \Big | ~ d \in \D
			~\text{s.t}~ \proj_x(d) \neq 0 \right \},
		\end{equation}
		is a PSS of $\T_x\M$ with $\cm_x(\fP_x(\D)) \ge \cm(\D)$.
	\end{proposition}}
\begin{proof}
	\rrev{Let $v\in \T_x\M\setminus \{0\}$. Since $D$ is a PSS, we know
		that
		\begin{equation*}
			\label{eq:projPSSTxMaux1}
			0
			< \cm(D)
			\le \sup_{d \in \D} \frac{\langle v, d \rangle_x}{ \|v\|_x\|d\|_x}
			= \sup_{\substack{d \in \D \\ \langle v,d \rangle_x > 0 }}
			\frac{\langle v, d \rangle_x}{ \|v\|_x\|d\|_x},
		\end{equation*}
		where the last inequality holds because the supremum on the left-hand
		side is positive. Using now that
		$\langle v,d \rangle_x = \langle v, \proj_x(d) \rangle_x$ and
		$\|d\|_x \geq \|\proj_x(d)\|_x$, we obtain
		\begin{equation*}
			\label{eq:projPSSTxMaux2}
			\sup_{\substack{d \in \D \\ \langle v,d \rangle_x > 0 }}
			\frac{\langle v,d \rangle_x}{ \|v\|_x\|d\|_x}
			=
			\sup_{\substack{d \in \D \\ \langle v,\proj_x(d) \rangle_x > 0 \\
					\|\proj_x(d)\|_x \neq 0}} \frac{\langle v, \proj_x(d) \rangle_x}{ \|v\|_x\|d\|_x}
			\le
			\sup_{\substack{d \in \D \\ \langle v,\proj_x(d) \rangle_x > 0 \\
					\|\proj_x(d)\|_x \neq 0}} \frac{\langle v, \proj_x(d) \rangle_x}{ \|v\|_x\|\proj_x(d)\|_x}.
		\end{equation*}
		Using~\eqref{eq:projPSSTxM}, it then follows that
		\begin{equation*}
			\sup_{\substack{d \in \D \\ \langle v,\proj_x(d) \rangle_x > 0 \\
					\|\proj_x(d)\|_x \neq 0}} \frac{\langle v, \proj_x(d) \rangle_x}{ \|v\|_x\|\proj_x(d)\|_x}
			\le
			\sup_{\substack{d \in \D \\ \|\proj_x(d)\|_x \neq 0 }}
			\frac{\langle v, \proj_x(d) \rangle_x}{ \|v\|_x\|\proj_x(d)\|_x}
			= \sup_{\bar d \in \fP_x(\D)} \frac{\langle v, \bar d \rangle_x}{ \|v\|_x\|\bar d\|_x}.
		\end{equation*}
		Overall, for any $v \in \T_x\M\setminus \{0\}$, we have shown that
		\begin{equation}
			\label{eq:cosine_inequality_proj}
			0
			<
			\cm(D)
			\le
			\sup_{d \in \D} \frac{\langle v, d \rangle_x}{ \|v\|_x\|d\|_x}
			\le
			\sup_{\bar d \in \fP_x(\D)} \frac{\langle v, \bar d \rangle_x}{ \|v\|_x\|\bar d\|_x}.
		\end{equation}
		Taking the infimum over $v \in \T_x\M \setminus \{0\}$
		in~\eqref{eq:cosine_inequality_proj} finally gives
		\begin{equation*}
			0
			<
			\cm(D)
			\le
			\inf_{ v \in \T_x\M \setminus \{0\}} \sup_{d \in \D}
			\frac{\langle v, d \rangle_x}{ \|v\|_x\|d\|_x}
			\le
			\inf_{v \in \T_x\M \setminus \{0\}} \sup_{\bar d \in \fP_x(\D)}
			\frac{\langle v, \bar d \rangle_x}{ \|v\|_x\|\bar d\|_x}
			=
			h\, \cm_x(\fP_x(\D)),
		\end{equation*}
		proving the desired result.
	}
\end{proof}

\rev{It follows from Proposition~\ref{prop:projPSSTxM} that the set
	$\fP_x(\D)$ satisfies Assumption~\ref{assu:pss} (with $\Dmin=\Dmax=1$) when
	$\cm(\D) \ge \kappa$. Building on this proposition, we now propose a PSS
	construction for the manifold $\M$.}

\begin{definition}
	\label{def:projPSSM}
	\rrev{Let $\D$ be a PSS of the ambient space $\R^n$.} Consider the family
	\begin{equation}
		\label{eq:projPSSM}
		\PM(D) := (\fP_x(\D))_{x \in \M},
	\end{equation}
	where $\fP_x(\D)$ is given by~\eqref{eq:projPSSTxM} for every $x \in \M$.
	The family~\eqref{eq:projPSSM} defines a PSS of $\M$ called \emph{projected
		PSS} induced by $\D$.
\end{definition}

Proposition~\ref{prop:projPSSTxM} implies that any PSS of $\M$ built
through~\eqref{eq:projPSSM} from a PSS $\D$ satisfies
$\cm_\M \left( \PM(\D) \right)\geq \cm(\D)$.
For instance, using the PSS $\pssa$ defined in $\R^n$,
one can build a family $\pssaproj{}{}$ that is a PSS of $\M$ with
$\cm_\M \left(\pssaproj{}{}\right) \ge \tfrac{1}{\sqrt{n}}$.

\begin{remark}
	\label{rem:VScmSubspace}
	The cosine measure of a projected PSS $\fP_x(\D)$ \rev{may differ} from the cosine measure of $\D$
	relative to the tangent space $\T_x\M$ recently proposed by Audet et
	al.~\cite{CAudet_WHare_GJarryBolduc_2025}. Indeed, all directions in $\fP_x(\D)$ are
	normalized, thus our cosine measure~\eqref{eq:cm_TxM} of $\fP_x(\D)$ only
	depends on the directions of the projections of the vectors in $\D$. On the contrary,
	the cosine measure of $\D$ relative to the tangent space $\T_x\M$ would depend on those
	norms. \rev{Normalizing projections in~\eqref{eq:projPSSTxM} ensures that even directions
		nearly orthogonal to the tangent space can significantly contribute to positively spanning
		the tangent space.}
\end{remark}

Projected PSSs were used in \rrev{Kungurtsev} et al.~\cite{VKungurtsev_FRinaldi_DZeffiro_2024}
in their numerical experiments, for which the ambient space could always be viewed as
$\R^n$ for some $n$, and projection was straightforward or provided through existing
manifold optimization packages. Note, however, that projections from the ambient space are
not always easy to compute. In addition, using projected PSSs generally leads to sets with
more vectors than in the intrinsic setting, without guarantees that the resulting set
will possess a good cosine measure in $\T_x\M$.

\begin{table}[htb!]
	\centering
	\caption{Illustration of an intrinsic PSS built from the PSS $\pssa$ of $\R^m$
		and a projected PSS built from the PSS $\pssa$ of $\R^n$.}

	\label{table:PSSa_intr_proj}

	\begin{tabular}{{lccc}} 
		\toprule
		PSS symbol                  & $\pssaintr{\B_x}$    & $\pssaproj{x}{}$         \\
		\midrule
		representation  in $\T_x\M$ &
		\begin{tikzpicture}[scale = 0.6, baseline={(current bounding box.center)}]

			\begin{scope}[xshift = 0,yshift = 0, scale = 1.35]
				\draw[gray, fill=gray!20, opacity=0.5] (-1.0, -1.0) -- (-1.0, 1.0) -- (1.0, 1.0) -- (1.0, -1.0) -- cycle;
				\node at (1.0,1.0) {$\T_x\M$};
			\end{scope}

			\draw[->,thick] (0,0) -- (0,1);
			\draw[->,thick] (0,0) -- (0,-1);
			\draw[->,thick] (0,0) -- (1,0);
			\draw[->,thick] (0,0) -- (-1,0);
		\end{tikzpicture}
		                            &
		\begin{tikzpicture}[scale = 0.6, baseline={(current bounding box.center)}]
			\begin{scope}[xshift = 0,yshift = 0, scale = 1.35]
				\draw[gray, fill=gray!20, opacity=0.5] (-1.0, -1.0) -- (-1.0, 1.0) -- (1.0, 1.0) -- (1.0, -1.0) -- cycle;
				\node at (1.0,1.0) {\rev{$\T_x\M$}};
			\end{scope}
			\draw[->,thick] (0,0) -- (1,0);
			\draw[->,thick] (0,0) -- (1/6,0.986013297);
			\draw[->,thick] (0,0) -- (-1/6,0.986013297);
			\draw[->,thick] (0,0) -- (-1,0);
			\draw[->,thick] (0,0) -- (1/6,-0.986013297);
			\draw[->,thick] (0,0) -- (-1/6,-0.986013297);

		\end{tikzpicture}               \\
		\midrule
		cardinality                 & $2m$                 & $\le 2n$                 \\
		\midrule
		cosine measure $\cmx$       & $\frac{1}{\sqrt{m}}$ & $\ge \frac{1}{\sqrt{n}}$ \\
		\bottomrule
	\end{tabular}
\end{table}

\begin{figure}[htb!]
	\centering
	\begin{tikzpicture}
		\begin{scope}[scale = 0.8]
			\tikzstyle{point}=[circle,fill=black,inner sep=1.5pt]

			\coordinate (o) at (0,0);
			\begin{scope}[xshift = 30,yshift = -5]
				\coordinate (a) at (-3.5,-3.) {};
				\coordinate (b) at (-2.5,1.) {};
				\coordinate (c) at (3,2.5) {}; 
				\coordinate (d) at (2,-1.2) {};
				\coordinate (cbis) at ($(c)!1/2!(d)!0.3!(o)$){};
				\coordinate (dbis) at ($(d)!1/2!(a)!0.5!(o)$){};
			\end{scope}
			\draw[-,thick] [smooth , tension = 0.5]
			plot coordinates{(a) (b) (c) };

			\node at ($(a)-(0.1,0.2)$) {$\mathcal{M}$};
			\draw[-, thick] (d) .. controls (dbis) .. (a)  node[right] {};
			\draw[-, thick] (c) .. controls (cbis) .. (d)  node[right] {};

			\draw[gray, fill=gray!20, opacity=0.5] (-1.5, 0.8) -- (2.5, 1.5) -- (1.5, -1.) -- (-2.5, -1.7) -- cycle;
			\node at (2.3, 1.6) {$T_{x} \mathcal{M}$};

			\node[point] (o) at (0,0) {};
			\node[left] at (-0.1,0.15) {$x$};

			\draw[->, thick] (o) -- (1.2, 0.6) node[above right] {};

			\draw[->, thick] (o) -- (-0.3, 0.7);
			\draw[->, thick] (o) -- (-0.1, 0.8) ;

			\draw[->, thick](o) -- (-1.2, -0.6);

			\draw[->, thick] (o) -- (0.4, -0.9);

			\draw[circle,blue, opacity = 0.7, thick] (0.18,-0.9) circle[x radius = 0.4, y radius = 0.15];
			\draw[-, blue] (0.2,-1.1) -- (0.7, -2.7) node[right,scale = 0.9] {$\substack{\text{Redundant} \\ \text{directions}}$};

			\draw[->, thick] (o) -- (0.1, -1.0);


			\begin{scope}[xshift = -180,yshift = 0, scale = 1.2]
				\draw[->, thick] (0,0) -- (0,-1) ;
				\draw[->, thick] (0,0) -- (1, 0) ;
				\draw[->, thick] (0,0) -- (-1, 0);
				\draw[->, thick] (0,0) -- (0,1);
				\draw[->, thick] (0,0) -- (0,-1);
				\draw[->, thick] (0,0) -- (1, 0);
				\draw[->, thick] (0,0) -- (0.78,0.5);
				\draw[->, thick] (0,0) -- (-0.78,-0.5);
				\node at (0.7, -0.5) {\rev{$\pssa$}};
			\end{scope}

			\begin{scope}[xshift = -120,yshift = 0, scale = 1.2]
				\draw[->, very thick, black] (-0.5,0.).. controls (0.9,0.5) .. (2.5,0) ;
				\node[fill=white, opacity=0, text opacity=1,scale = 1] at (1,0.6) {\textcolor{black}{$\mathrm{proj}_x$}};
			\end{scope}
		\end{scope}
	\end{tikzpicture}
	\caption{Projection of a positive spanning set $\pssa$ in $\R^3$ to a
		tangent space $\T_x\M$. The projection produces close or redundant directions,
		that affect the cardinality but may not increase the cosine measure
		significantly.}
	\label{fig:projecting_PSSs}
\end{figure}

Table~\ref{table:PSSa_intr_proj} and Figure~\ref{fig:projecting_PSSs}
illustrate possible drawbacks \rev{of} using projected PSSs compared to
intrinsic PSSs. On one hand, the number of vectors in the projected PSS may be
as large as $2n$, thus possibly much larger than $2m$. On the other hand, the
projected PSS may not have significantly better geometry than an intrinsic
PSS, because some projected directions may be close to one another. In the
next section, we provide a formal justification for this illustration in the
case of the sphere manifold.

\section{Polling directions on the sphere manifold}
\label{sec:study_cm_sphere}

In this section, we focus on building PSSs for the sphere manifold
$\mathbb{S}^{n-1}= \{x \in \R^n\ |\ \| x \|_2^2 = 1 \}$. This manifold is an
embedded submanifold of $\R^n$ of dimension $n-1$. For every
$x \in \mathbb{S}^{n-1}$, the tangent space at $x$ is given by
\begin{equation}
	\label{eq:tangent_sphere}
	\T_x\mathbb{S}^{n-1} = \left \{ u \in \R^n\ \middle|\ x^\T u = 0 \right \},
\end{equation}
and the corresponding projection operator is the mapping
$\proj_x : y  \to y - x x^\T y$. We will consider PSSs for $\mathbb{S}^{n-1}$
built from $\pssa$, and show that the intrinsic approach has better
complexity measure than its projected counterpart.

Section~\ref{ssec:theory_cosine_measure_sphere} is dedicated to estimating
the cosine measure and the complexity measure of the PSSs at hand for a
particular example. Experiments reported in Section~\ref{ssec:sphere:num}
suggest that our theoretical findings hold in a more general setting.

\subsection{Cosine and complexity measures for a projected PSS}
\label{ssec:theory_cosine_measure_sphere}

Given the PSS $\pssa$ in $\R^n$, we wish to \rrev{estimate the} complexity measure
of the projected PSS $\PM(\pssa)$. To this end, we focus first on the cosine
measure of this set, which we reformulate as the maximization of a convex
function over a polytope.

\begin{lemma} \label{lemma:pssa_sphere_cm_reformulation}
	For any $x \in \bbS^{n-1}$ such that $|x_i| < 1$ for every $i=1,\dots,n$,
	define the norm $\|\cdot\|_{\infty,x}$ by
	\[
		\|u\|_{\infty,x} = \max_{i=1,\dots,n} \frac{|u_i|}{\sqrt{1-x_i^2}}
		\quad
		\forall u \in \R^n.
	\]
	Then, we have
	\begin{align}
		\cm_x(\pssaproj{x}{})
		 & = \left (
		\max_{u \in \R^n} \left \{
		\|u\|_2\ \mathrm{s.t.}\ \|u\|_{\infty, x} \le 1, u^\T x = 0
		\right \} \right )^{-1} \label{eq:cmprojPSS1_eq2}.
	\end{align}
\end{lemma}
\begin{proof}
	By definition of the cosine measure, we have
	\[
		\cm_x(\pssaproj{x}{})
		=
		\min_{\substack{u \in \R^n \\ \|u\|_2 = 1\\u^\T x = 0}}
		\max_{d \in \pssa} \frac{\langle \proj_x(d),u \rangle}{\| \proj_x(d)\|}.
	\]
	Recall that $\pssa$ consists of the coordinate vectors and their
	negatives in $\R^n$. For any vector $\pm e_i \in \pssa$, applying
	the projection gives:
	\begin{align*}
		 & \proj_x(\pm e_i) = \pm (e_i -x_i x) \neq 0,
		 & \|\proj_x(\pm e_i)\|_2^2 = (1-x_i^2)^2 + \sum_{j \neq i} x_i^2 x_j^2 = 1-x_i^2,
	\end{align*}
	where the last equality uses $x \in \mathbb{S}^{n-1}$.
	Thus,
	\[
		\cm_x(\pssaproj{x}{})
		= \min_{\substack{u \in \R^n \\ \|u\|_2 = 1\\u^\T x = 0}} \max_{1 \leq i \leq n}
		\max \left ( \frac{\langle e_i - x_i x,u \rangle}{\sqrt{1-x_i^2}} ,
		-\frac{\langle e_i - x_i x,u \rangle}{\sqrt{1-x_i^2}} \right )
		= \min_{\substack{u \in \R^n \\ \|u\|_2 = 1\\u^\T x = 0}} \max_{1 \leq i \leq n}
		\frac{|\langle e_i - x_i x,u \rangle|}{\sqrt{1-x_i^2}}.
	\]
	Using $u^\T x=0$, we have $\langle e_i-x_i x,u \rangle = u_i$ for every $i$.
	Therefore, we obtain
	\begin{equation}
		\label{eq:cmprojPSS1_eq1}
		\cm_x(\pssaproj{x}{})
		=
		\min_{\substack{u \in \R^n \\ \|u\|_2 = 1\\u^\T x = 0}} \|u\|_{\infty,x}.
	\end{equation}

	We now proceed to proving~\eqref{eq:cmprojPSS1_eq2}
	from~\eqref{eq:cmprojPSS1_eq1}. Using homogeneity of the \rev{norm
	$\| \cdot\|_{\infty, x}$} gives
	\begin{equation}
		\label{eq:homogeneousmin}
		\min_{\substack{u \in \R^n \\ \|u\|_2 = 1\\u^\T x = 0}}
		\|u\|_{\infty,x}
		=
		\min_{\substack{u \in \R^n \setminus \{ 0 \} \\ u^\T x = 0}}
		\left \| \frac{u}{\|u\|_2} \right \|_{\infty,x}
		=
		\min_{\substack{u \in \R^n \setminus \{ 0 \} \\ u^\T x = 0}}
		\frac{\|u\|_{\infty,x}}{\|u\|_2}\rev{.}
	\end{equation}
	Noticing that for any $u \in \T_x \bbS^{n-1}$, we have
	\begin{equation*}
		\min_i (1-x_i^2)^{-\frac 12} \|u\|_{\infty}
		\le
		\|u\|_{\infty,x}
		\le
		\max_i (1-x_i^2)^{-\frac 12} \|u\|_{\infty},
	\end{equation*}
	and using $0<\|u\|_{\infty} \le \|u\|_2 \le \sqrt{n}\|u\|_{\infty}$,
	we arrive at
	\begin{equation}
		\label{eq:ratio_norms_bounded}
		\frac{1}{\sqrt{n}}\min_i (1-x_i^2)^{-\frac 12}
		\le
		\frac{\|u\|_{\infty,x}}{\|u\|_2}
		\le
		\max_i (1-x_i^2)^{-\frac 12}.
	\end{equation}
	Equation~\eqref{eq:ratio_norms_bounded} guarantees that the quantity
	\begin{equation}
		\label{eq:inverse ratio}
		\max_{\substack{u \in \R^n \setminus \{0\} \\ u^\T x = 0}}
		\frac{\|u\|_2}{\|u\|_{\infty,x}}
	\end{equation}
	is finite and positive. Moreover, we have
	\begin{equation}
		\label{eq:inversepropty}
		\min_{\substack{u \in \R^n \setminus \{ 0 \} \\ u^\T x = 0}}
		\frac{\|u\|_{\infty,x}}{\|u\|_2}
		=
		\left[\max_{\substack{u \in \R^n \setminus \{0\} \\ u^\T x = 0}}
		\frac{\|u\|_2}{\|u\|_{\infty,x}}\right]^{-1}
	\end{equation}
	since both objective functions are positive on the set
	$\{u \in \R^n \setminus \{0\} ~|~ u^\T x = 0\}$. Finally, we
	can reformulate~\eqref{eq:inverse ratio} using homogeneity as
	\begin{equation}
		\label{eq:homogeneousmax}
		\max_{\substack{u \in \R^n \setminus \{ 0 \} \\ u^\T x = 0}}
		\frac{\|u\|_2}{\|u\|_{\infty,x}}
		=
		\max_{\substack{u \in \R^n \setminus \{ 0 \} \\ u^\T x = 0}}
		\left \| \frac{u}{\|u\|_{\infty,x}} \right \|_2
		=
		\max_{\substack{u \in \R^n \\ \|u\|_{\infty,x}= 1\\ u^\T x = 0}}
		\|u\|_2.
	\end{equation}
	Combining~\eqref{eq:cmprojPSS1_eq1}, \eqref{eq:homogeneousmin},
	\eqref{eq:inversepropty} and~\eqref{eq:homogeneousmax}, we arrive at the
	desired conclusion.
\end{proof}

Using Lemma~\ref{lemma:pssa_sphere_cm_reformulation}, we can provide a bound
on the cosine measure of $\pssaproj{x}{}$ for any $x$, and even instantiate it
for specific choices of the point $x$.
\begin{theorem}
	\label{thm:bound_projcm_pssa}
	For any $x \in \bbS^{n-1}$, we have
	\begin{equation}
		\label{eq:bound_cm_pssaproj_sphere}
		\frac{1}{\sqrt{n-1}}
		\le \cm_x(\pssaproj{x}{}) \le
		\frac{1}{\sqrt{n-2+\|x\|_\infty^2}}
		\le \frac{\sqrt{n}}{n-1}.
	\end{equation}

	In addition, if $x$ has $k$ components equal to $\pm \tfrac{1}{\sqrt{k}}$
	for some $k \in \{1,\dots,n\}$,
	\begin{equation}
		\label{eq:cm_pssaproj_partpoints}
		\cm_x(\pssaproj{x}{}) =
		\begin{cases}
			\frac{1}{\sqrt{n-2 + \frac{1}{k}}} & \text{if } k\text{ is odd},
			\\
			\frac{1}{\sqrt{n-1}}               & \text{if } k \text{ is even}.
		\end{cases}
	\end{equation}
\end{theorem}
\begin{proof}
	Let $x \in \bbS^{n-1}$. Without loss of generality, we prove the result
	assuming that $x_1 \ge x_2 \ge \cdots \ge x_n \ge 0$, and we let $k$ be the
	first index in $\{1,\dots,n\}$ such that $x_{k+1}=0$, with $k=n$ if all
	coordinates of $x$ are positive.

	If $k=1$, then $x=e_1$, in which case $\pssaproj{x}{} = \pssaintr{B}$, where
	$B = (e_2, \dots, e_n)$. If $k>1$, we define
	\begin{align}
		 & \calQ_x = \left\{ u \in \R^n\ \middle|\ \|u\|_{\infty,x} \le 1\ \mathrm{and}\ u^\T x = 0 \right\} \\
		 & \tau(x):= \max_{u \in \calQ_x} \|u\|_2. \label{eq:2norm_maximization_cmPSS1sphere}
	\end{align}
	From Lemma~\ref{lemma:pssa_sphere_cm_reformulation}, we know that
	$\cm_x(\pssaproj{x}{})=\tau(x)^{-1}$. It thus suffices to bound $\tau(x)$
	to obtain bounds on $\cm_x(\pssaproj{x}{})$.

	The value $\tau(x)$ is the maximum of a convex function over the polytope $\calQ_x$,
	hence it is attained at extreme points $\calE(\calQ_x)$ of
	$\calQ_x$~\cite[Theorem 32.3]{RTRockafellar_1970}. For any $v \in \calE(\calQ_x)$,
	we have
	\begin{equation}
		\label{eq:propextremeQx}
		|v_i| = 1\ \forall i \ge k+1
		\quad \mbox{and} \quad
		|v_i| \neq \sqrt{1-x_i^2}\ \mbox{for at most \rev{one} } i \in \llbracket 1,k \rrbracket.
	\end{equation}
	Indeed, suppose that there exists $i \ge k+1$ such that $|v_i| \neq 1$. Let
	$v^+ = v + \eta e_i$ and $v^{-} = v - \eta e_i$ with $\eta = \frac{1-|v_i|}{2}$. Then
	both $v^+$ and $v^-$ belong to $\calQ_x$ while being distinct of $v$, hence
	$v = \frac{v^+ + v^-}{2}$, which is a contradiction.
	Similarly, suppose that there \rev{exist} two distinct indices
	$i,j \in \llbracket 1,k \rrbracket $ such that $|v_i| \neq \sqrt{1-x_i^2}$ and
	$|v_j| \neq \sqrt{1-x_j^2}$. Define
	\begin{equation*}
		v^+ = v + \frac{\eta}{x_i}e_i - \frac{\eta}{x_j}e_j, ~~~
		v^- = v - \frac{\eta}{x_i}e_i + \frac{\eta}{x_j}e_j
	\end{equation*}
	with $\eta = \frac{1}{2}\min(x_i(\sqrt{1-x_i^2}-|v_i|), x_j(\sqrt{1-x_j^2}-|v_j|))$,
	we again obtain that both $v^+$ and $v^-$ lie in $\calQ_x$ while being
	different from $v=(v^+ + v^-)/2$, which is \rev{another} contradiction.

	Using~\eqref{eq:propextremeQx}, we obtain that $v$ can be written as
	\begin{equation}
		\label{eq:extremeptstruct}
		v =
		\begin{pmatrix}
			\eps_1 \sqrt{1-x_1^2}
			 & \cdots
			 & \eps_{i-1} \sqrt{1-x_{i-1}^2}
			 & \nu
			 & \eps_{i+1} \sqrt{1-x_{i+1}^2}
			 & \cdots
			 & \eps_{n} \sqrt{1-x_{n}^2}
		\end{pmatrix}^\T ,
	\end{equation}
	for some $i \in \llbracket 1,k \rrbracket$, where $\eps_j \in \{-1,1\}$
	for any $j \neq i$ and
	$\nu = -\tfrac{1}{x_i}\sum_{j \neq i} \eps_j x_j \sqrt{1-x_j^2} \le \sqrt{1-x_i^2}$.
	We then reformulate~\eqref{eq:2norm_maximization_cmPSS1sphere} as a problem
	over the extreme points of $\calQ_x$:
	\begin{equation}
		\label{eq:reformtau}
		\tau(x) = \max_{v \in \calE(\calQ_x)} \|v\|_2
		\rrev{~= \max_{i \in \llbracket 1,k \rrbracket} \tau_i(x),}
	\end{equation}
	where
	\begin{equation}
		\label{eq:tauk}
		\begin{array}{rcl}
			\rrev{\tau_i(x)} = \max_{\substack{(\eps_j)_{j \neq i} \in \{-1,1\}^{k-1} \\ \nu \in \R}}
			\sqrt{(\sum_{j \neq i} 1-x_j^2) + \nu^2}
			 & \mbox{s.t.}
			 & \nu= - \frac{1}{x_i} \sum_{\substack{j=1                               \\j \neq i}}^{k} \epsilon_j x_j \sqrt{1-x_j^2},
			\\
			 &
			 & |\nu| \leq \sqrt{1-x_i^2}.
		\end{array}
	\end{equation}
	The first maximum in~\eqref{eq:reformtau} accounts for the choice of the non-active inequality in
	$\llbracket 1,k \rrbracket $ and the second in~\eqref{eq:tauk} corresponds to the maximization of the objective over extreme points having the structure~\eqref{eq:extremeptstruct} for
	that particular $i$. Squaring the objective and using $\|x\|_2^2 = 1$, we obtain
	\begin{equation} \label{eq:pm1_formulation_cm_projpss1sphere}
		\tau(x)^2 = n-2 + \max_{i \in \llbracket 1,k \rrbracket  } [x_i^2 + \alpha_i^2],
	\end{equation}
	where
	\[
		\alpha_i^2 =
		\max_{(\epsilon_j)_{j \neq i} \in \{-1,1\}^{k-1}, \alpha \in \R} \left\{ \alpha^2
		\quad \mbox{s.t.} \quad
		\alpha = - \frac{1}{x_i} \sum_{\substack{j=1 \\j \neq i}}^{k} \epsilon_j x_j \sqrt{1-x_j^2},
		|\alpha| \leq \sqrt{1-x_i^2}
		\right\}.
	\]
	Using the bound $\alpha_i^2 \le 1-x_i^2$
	in~\eqref{eq:pm1_formulation_cm_projpss1sphere} yields
	\begin{equation}
		\label{eq:boundtaux}
		n-2 + \|x\|_\infty^2
		\le
		\tau(x)^2
		\le
		n-2+ \max_{1\leq i\leq k}x_i^2 +(1-x_i^2) = n-1,
	\end{equation}
	which leads to
	\[
		\frac{1}{\sqrt{n-1}}
		\le
		\cm_x(\pssaproj{x}{})
		\le
		\frac{1}{\sqrt{n-2+\|x\|_\infty^2 }}.
	\]
	It suffices to note that
	$\|x\|_\infty \geq \frac{1}{\sqrt{n}}\|x\|_2=\frac{1}{\sqrt{n}}$
	to obtain~\eqref{eq:bound_cm_pssaproj_sphere}.

	Suppose now that $x$ has $k$ components equal to $\pm \tfrac{1}{\sqrt{k}}$. To
	prove~\eqref{eq:cm_pssaproj_partpoints}, we consider two cases. If $k \ge 2$ is odd, then
	$\llbracket 1, k \rrbracket \setminus \{i\}$ has even cardinality. It follows
	that optimality for~\eqref{eq:pm1_formulation_cm_projpss1sphere} is
	attained by having an equal number of $\epsilon_j$ equal to $+1$ and $-1$,
	so as to respect the bound on $|\alpha|$. Thus, $\alpha_i = 0$ for every $i$ and
	$\tau(x)^2 = n-2+\|x\|_\infty^2$, proving~\eqref{eq:cm_pssaproj_partpoints}
	in that case.
	If $k\geq 2$ is even, then $\llbracket 1, k \rrbracket \setminus \{i\}$ has odd
	cardinality. Consequently, the bound on $|\alpha|$ can only be satisfied when
	the numbers of positive and negative values in $\{\epsilon_j\}_{j \neq i}$ differ by one.
	As a result, we must have $\alpha^2_i =1-x_i^2$ for every $i$ and $\tau(x)^2 = n-1$,
	which leads to the desired conclusion.
\end{proof}

\rev{A direct} consequence of Theorem~\ref{thm:bound_projcm_pssa} is that the PSS
$\pssaproj{}{}$ of $\bbS^{n-1}$ satisfies
\begin{equation}
	\label{eq:boundcmpssaproj}
	\frac{1}{\sqrt{n-1}}
	\le \cm_{\bbS^{n-1}}(\pssaproj{}{}) \le
	\frac{\sqrt{n}}{n-1}.
\end{equation}
We can now compare the projected PSS induced on $\bbS^{n-1}$ by $\pssa$ to
an intrinsic PSS induced by $\pssa$ in terms of complexity measure.
\begin{corollary}
	\label{coro:compintrprojsphere}
	Let $\BM = (\B_x)_{x \in \bbS^{n-1}}$ be a family of \rrev{orthonormal} bases
	adapted to the manifold $\bbS^{n-1}$. Let $\pssaintr{\BM}$ be the intrinsic
	PSS induced by $\BM$ and $\pssa$ (as a PSS in $\R^{n-1}$), and
	$\pssaproj{}{}$ be the projected PSS of $\bbS^{n-1}$ induced by $\pssa$ (as
	a PSS in $\R^n$). Then,
	\begin{equation}
		\label{eq:compintrprojsphere}
		\chi_{\bbS^{n-1}}(\pssaproj{}{}) = 2n(n-1)
		\ge
		2(n-1)^2 =\chi_{\bbS^{n-1}}(\pssaintr{\BM}).
	\end{equation}
\end{corollary}
%
%
%
\begin{proof}
	\rev{
		We begin by proving the first equality in~\eqref{eq:compintrprojsphere}.
		By definition of a projected PSS, we have $|\pssaproj{x}{}| \le |\pssa| = 2n$
		for every $x \in \M$, while Theorem~\ref{thm:bound_projcm_pssa} gives\\
		$
			\sup_{x \in \bbS^{n-1}} (\cmx(\pssaproj{x}{}))^{-2} = n-1.
		$
		Therefore,
		\begin{equation*}
			\chi_{\bbS^{n-1}}(\pssaproj{}{}) \le 2n(n-1).
		\end{equation*}
		Consider now the special case of $x$ being the vector}
	\rrev{\begin{equation*}
			x:= \begin{cases}
				\frac{1}{\sqrt{n}}(1,\dots,1,1)   & \text{if $n$ is even,} \\
				\frac{1}{\sqrt{n-1}}(1,\dots,1,0) & \text{otherwise}.
			\end{cases}
		\end{equation*}}
	\rev{For all $i  \in \{1, \dots, \}$,
		$\|\proj_x(\pm e_i)\|_x \neq 0$, and as a result $|\pssaproj{x}{}| = 2n$,
		while $(\cmx(\pssaproj{x}{}))^{-2}=n-1$ per
		Theorem~\ref{thm:bound_projcm_pssa}. It follows that
		$\chi_{\bbS^{n-1}}(\pssaproj{}{}) = 2n(n-1)$.

		We now turn to the second equality. Consider a family of isomorphisms
		$(\isomintr{\B_x})_{x \in \M}$ defining the PSS $\pssaintr{\BM}$ as in
		Definition~\ref{def:intrinsic_PSS}. By definition, for any $x \in \M$, we
		have
		\[
			\cm_x(\pssaintr{\B_x})
			= \cm_x(\isomintr{\B_x}(\pssa))
			= \cm(\pssa)
			=1/\sqrt{n-1},
		\]
		where the last $\pssa$ refers to the corresponding PSS in $\R^{n-1}$.
		Similarly, we find that
		$|\pssaintr{\B_x}| = |\isomintr{\B_x}(\pssa)| = |\pssa| = 2(n-1)$. We
		thus arrive at
		\[
			\chi_{\bbS^{n-1}}(\pssaintr{\BM})
			=\sup_{x \in \M}|\pssaintr{\B_x}|
			\left(\inf_{x \in \M}\cm_x(\pssaintr{\B_x})\right)^{-2}
			=2(n-1)^2,
		\]
		proving the desired result.}
\end{proof}

Corollary~\ref{coro:compintrprojsphere} shows that the intrinsic PSS
approach yields a better complexity measure than the projected PSS
approach. Although encouraging, our theory \rev{does} not extend easily to other PSSs
on the sphere, as they heavily rely on properties of $\pssa$. This further
supports the use of intrinsic PSSs as they readily come with complexity measure
guarantees.

\subsection{Distribution of cosine measure values for a projected PSS}
\label{ssec:sphere:num}

In this section, we numerically evaluate the cosine measure for the projected
PSS $\pssaproj{}{}$. Although the results of
Section~\ref{ssec:theory_cosine_measure_sphere} \rev{fix} a worst case for this
PSS, we aim at illustrating the distribution of cosine measures of
$\pssaproj{x}{}$ across tangent spaces.

To this end, we first compute the value of $\cmx(\pssaproj{x}{})$ for $n=3$
through a discretization of the sphere via spherical coordinates. After
computing $\pssaproj{x}{}$ for each $x$ in our discretization, we decompose
the set using an \rrev{orthonormal} basis of $\T_x\bbS^{2}$, then apply the
cosine measure algorithm of Hare and Jarry-Bolduc~\cite{WHare_GJarryBolduc_2020}
to the resulting matrix.

\begin{figure}[h!]
	\centering
	\includegraphics[scale = 0.35]{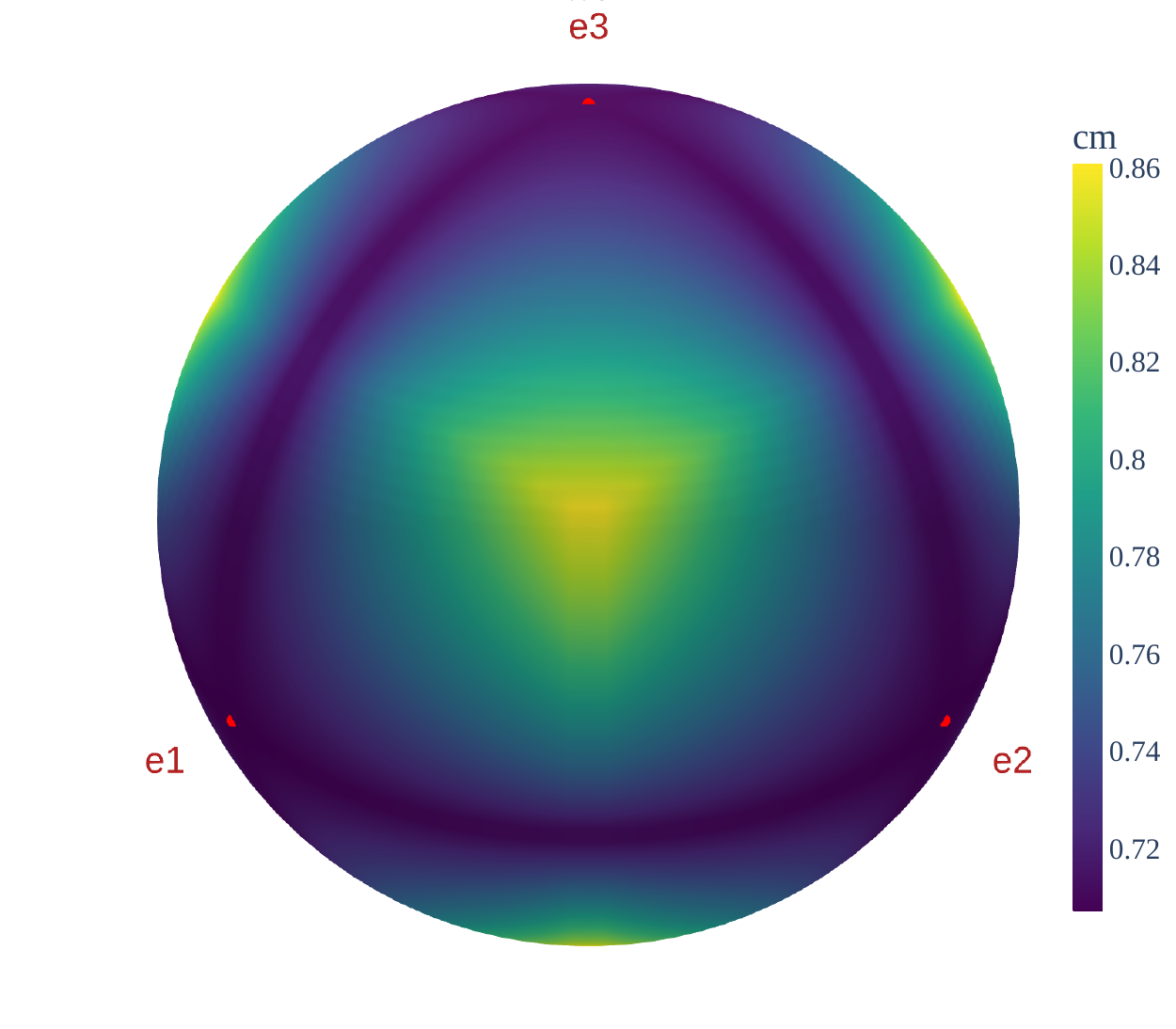}
	\includegraphics[scale = 0.35]{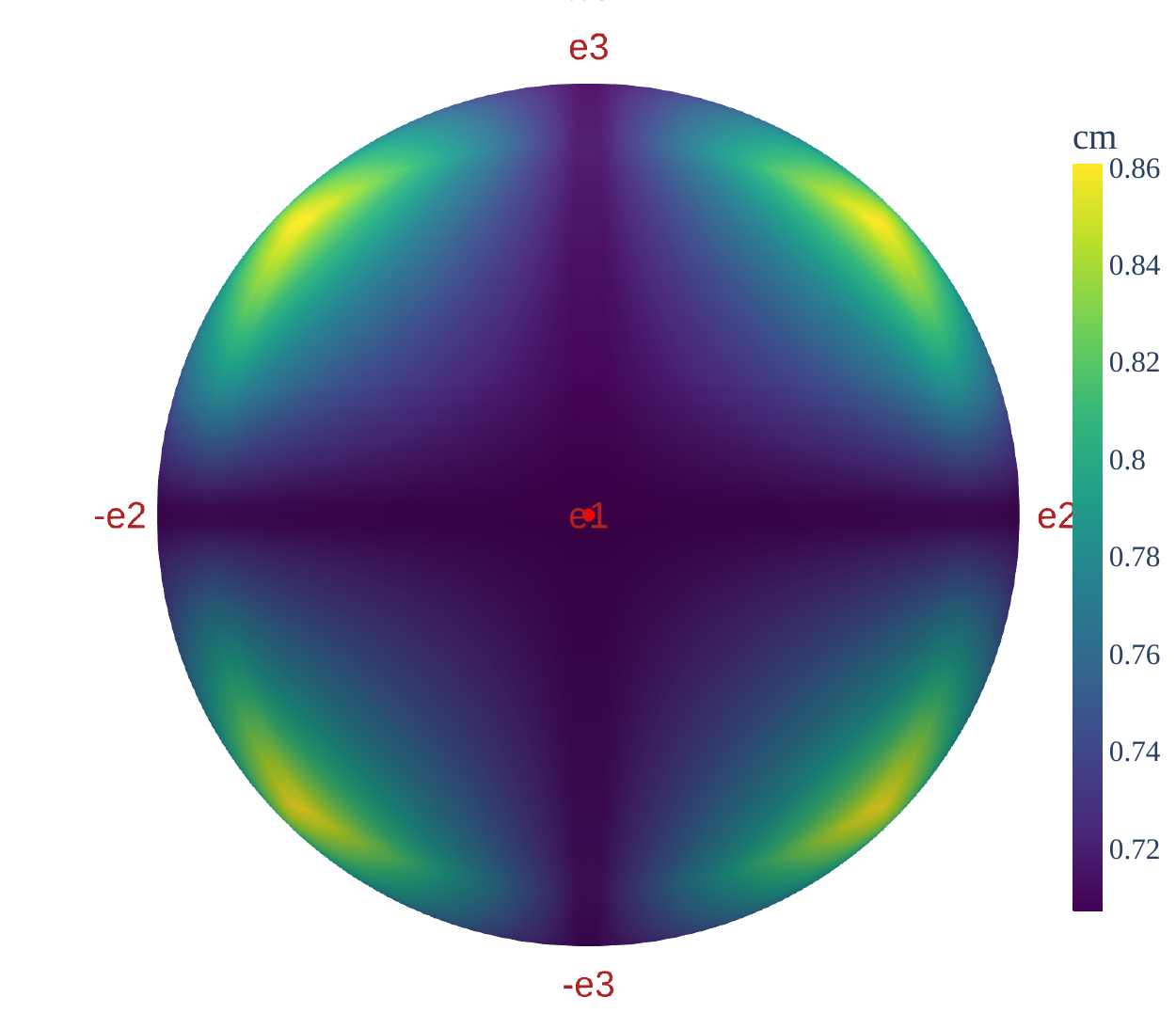}
	\caption{\prev{Cosine measure of the projected PSS $\pssaproj{x}{}$ at various points on the sphere $\bbS^{2} \subset \R^3$, shown from two viewpoints.}}
	\label{fig:num:pssproj_cmheatmap}
\end{figure}

Figure~\ref{fig:num:pssproj_cmheatmap} presents a heatmap of this function.
The symmetries in Figure \ref{fig:num:pssproj_cmheatmap} confirm that the value of
$\cmx(\pssaproj{x}{})$ is \rev{independent of the sign} of the components in $x$.
We also remark that the specific vectors highlighted in
Theorem~\ref{thm:bound_projcm_pssa} represent accurately the range of
cosine measure values across the sphere.

Although the sphere $\bbS^{n-1}$ has no natural representation for $n > 3$, we are
still able to compute the cosine measure at selected points through the same process
\rev{as that described} above. Note that we use a QR decomposition of a random matrix
with $x$ as the first column to compute a basis of $\T_x\bbS^{n-1}$.
\rev{
	\begin{figure}[h!]
		\centering
		\includegraphics[scale = 1]{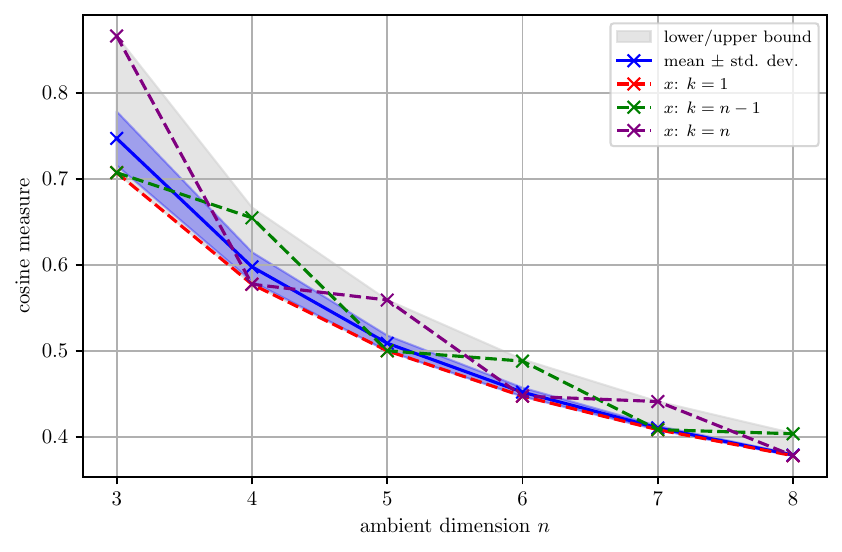}
		\caption{\prev{Cosine measure of the projected PSS $\pssaproj{x}{}$ at various} \rrev{points} \prev{on the sphere $\bbS^{n-1} \subset \R^n.$}} 
		\label{fig:num:cmrange_plot}
	\end{figure}}

Figure~\ref{fig:num:cmrange_plot} displays values of $\cmx(\pssaproj{x}{})$ at
selected points \prev{--} specifically those defined in Theorem \ref{thm:bound_projcm_pssa}
with $k=1,~n-1$ or $n$ for various values of $n$. As expected, those lie between the lower and
upper bounds provided by Theorem~\ref{thm:bound_projcm_pssa}. We also represent the mean
value of $\cmx(\pssaproj{x}{})$ computed over $100$ random points on the sphere \rrev{(generated
	using the \texttt{random\_point} method from the \texttt{pymanopt}
	library~\cite{JTownsend_NKoep_SWeichwald_2016})}. We observe that the mean shrinks towards the
lower bound as $n$ increases, while the gap between the lower and upper bound shrinks less
rapidly. Since the \rrev{lower bound} corresponds to the cosine measure of the intrinsic PSS, our
results suggest that the gap between projected and intrinsic PSSs built from $\pssa$ is
likely to have a practical impact. The next section aims at validating this hypothesis.

\section{Direct-search experiments}
\label{sec:num}

\rev{
	In this section, we compare direct-search methods based on intrinsic and
	projected PSSs, as defined in Sections~\ref{ssec:intrinsic_PSS}
	and~\ref{ssec:projected_PSS}. Our experiments aim not only at comparing those
	variants when induced by the same PSS structure, but also at assessing the
	effect of the codimension $n-m$ and the dimension $m$ in the difference in
	performance. Section~\ref{ssec:experimental_setup} presents our numerical
	study on toy manifold problems, and Section~\ref{ssec:num:obs} provides the
	results. Section~\ref{ssec:num:rotsync} then compares both direct-search
	variants on the problem of synchronization of rotations.
}

\subsection{Experimental setup}
\label{ssec:experimental_setup}

We conduct experiments using two kinds of manifold constraints. We first
consider product manifolds of the form
\begin{equation}
	\label{eq:num:embedsphere}
	\Msphere := \bbS^{(m+1)-1} \times \rrev {\left\{ 0_{\R^{n-m-1}} \right\},}
\end{equation}
that represent spheres of dimension $m$ embedded in $\R^n$. Our study
from Section~\ref{sec:dsdirs} considered solely the case $m=n-1$, where
the codimension and ambient dimension only differ by $1$. Using
manifolds~\eqref{eq:num:embedsphere} allows to study how the performance of
our algorithms changes with $n$ and $m$.
We also consider problems over subspace manifolds
\begin{equation}
	\label{eq:num:sev}
	\Mlinear \; := \; \left\{ x \in \R^n\ |\ x \in \spann(z_1,\dots,z_m)\right\},
\end{equation}
where $z_1,\dots,z_m$ are orthogonal vectors drawn randomly in $\R^n$. \prev{These manifolds} are fairly simple, allowing the dimension and
codimension to take almost any value. Note that the codimension for the
manifold~\eqref{eq:num:embedsphere} is at least $1$, whereas that for the
manifolds~\eqref{eq:num:sev} can be $0$. Nevertheless, we will use $n$ and $m$
indifferently for both manifolds, since our goal is to study variations of
the performance as $n-m$ (or $m$) increases.

Our test problems are inspired from barycenter and eigenvalue
problems~\cite{NBoumal_2023,XPennec_2006} that were previously used for
Riemannian direct search~\cite{VKungurtsev_FRinaldi_DZeffiro_2024}.
The first three problems are constrained to lie on the manifold $\Mlinear$.
For the first two problems, we draw $10$ points $(x_1, \dots, x_{10})$ at random
and consider the objective function
\begin{equation*}
	\sum_{i= 1}^{10} \|x_i-x\|_2^2,
\end{equation*}
whose minimizer is $x = \proj_\Mlinear\left (\frac{1}{10}\sum x_i \right )$.
In the first case, \rev{the points $x_1,\dots,x_{10}$ are drawn} in $\R^n$,
whereas in the second case these points are drawn in the linear subspace
$\Mlinear$. The third objective is a random $\tfrac{1}{10}$-strongly convex
quadratic. The fourth problem corresponds to minimizing the Rayleigh quotient
of a random symmetric matrix over $\Msphere$. Table~\ref{table:numsetup_problems}
summarizes our test problems and the associated manifolds. We generate $100$
instances of each problem using manifold dimensions $m \in \{2,4,8,16,32\}$
and codimensions $n-m$ ranging in \rrev{$\{0,2,4,8,16,32\}$ for $F$ and
	$\{1,2,4,8,16,32\}$ for $\Msphere$).}
By using this simple benchmark, we hope that differences in performance
will be essentially due to the polling strategies rather than the difficulty of a
given problem.
\begin{table}[htb!]
	\centering
	\caption{Benchmarking problems for direct-search}
	\label{table:numsetup_problems}
	\begin{tabular}{ccc}
		\toprule
		Manifold                                                                  & $f(x)$                                                   & Parameters                                     \\
		\midrule
		$\Mlinear = \spann(z_1, \dots, z_m)$                                      & $\sum_{i= 1}^{10} \|x_i-x\|_2^2$                         & $(x_i) \in (\R^n)^{10}$                        \\

		                                                                          & $\sum_{i= 1}^{10} \|x_i-x\|_2^2$                         & $(x_i) \in \Mlinear^{10}$                      \\

		                                                                          & $\frac{1}{2} \langle Ax,x \rangle - \langle b,x \rangle$ & $A \succcurlyeq \frac{1}{10}I_n ,~ b \in \R^n$ \\
		\midrule
		$\rev{\Msphere =}~ \mathbb{S}^{(m+1)-1} \times \rrev{\{0_{\R^{n-m-1}}\}}$ & $\langle x, Ax \rangle$                                  & $A \in \mathcal S ym(n)$                       \\
		\bottomrule
	\end{tabular}
\end{table}

We implemented Algorithm~\ref{algo:RDS} using $\gammadec = \frac{1}{2}$,
$\gammainc = 2$, $\alpha_{\max} = 1$, $c = 1$ and $x_0$ as a random point on the
manifold. We used six different polling choices, corresponding to $\pssaproj{}{}$,
$\pssbproj{}{}$, $\psscproj{}{}$, $\pssaintr{\BM}$, $\pssbintr{\BM}$ and
$\psscintr{\BM}$ in Section~\ref{sec:dsdirs}.
To construct projected PSSs, we rely on the canonical basis \prev{of} $\R^n$.
For every point $x \in \Mlinear$, we use $\B_x = (z_1,\dots,z_m)$, where
$(z_1,\dots,z_m)$ spans $F$, to build intrinsic PSSs at $x$. For every
$x \in \Msphere$, akin to Section~\ref{ssec:sphere:num}, we first build a
basis of $\bbS^{m}$ from the QR decomposition of a random $(m+1) \times (m+1)$
matrix where the first column is replaced by the first $m+1$ coefficients of $x$,
then augment the vectors in this basis with zero coordinates to form a basis
$\B_x$ of $\T_x\Msphere$.

We run all variants of Algorithm~\ref{algo:RDS} for a budget of
\rev{$100(m+1)$	function evaluations} for a problem with manifold dimension
$m$. \rev{This budget corresponds to $100$ simplex gradient estimates in a
	space of dimension $m$. A simplex gradient is an approximate gradient
	construction, and its basic cost ($m+1$ evaluations in dimension $m$)
	is used to normalize evaluation budget according to the
	dimension~\cite{JJMore_SMWild_2009}. Our results consider the dimension
	of the manifold ($m$) rather than that of the ambient space ($n \ge m$) which
	would result in a higher budget per problem. However, experiments conducted
	with the larger budget of $100(n+1)$ evaluations lead to the same conclusions.}

To compare our polling strategies, we make use of data
profiles~\cite{JJMore_SMWild_2009}. \rrev{We consider a set of problems $\mathcal{P}$
with associated objective functions $\{f_p\}_{p \in \mathcal{P}}$. Given a
problem $p$, a solver $s$, \rev{a tolerance $\tau>0$,} and a
proxy $f_p^*$ for the optimal value of this problem}, we let
$t_{p,s}^\tau$  denote the number of evaluations needed by solver $s$ to obtain
an iterate $x_k$ such that
\rrev{\[
		\frac{f_p(x_k) - f_p^*}{f_p(x_0)-f_p^*} \leq \tau.
	\]}
The data profile of solver $s$ is the cumulative distribution function
\[
	d_s(\alpha) = \frac{1}{|\mathcal{P}|} \left |
	\left \{ p \in \mathcal{P}: t_{p,s}^\tau \leq \alpha(m+1) \right \}\right |,
\]
counting the proportion of problems solved by solver $s$ for a given budget
$\alpha$. Higher curves \rev{correspond} to more efficient solvers.
\rev{The proxy $f_p^*$ is obtained by taking the best value between all
	solvers after running out of the budget. We plot all data profiles in
	Section~\ref{ssec:num:obs} with tolerance $\tau = 10^{-2}$, but our
	observations extend to other tolerances.}


Our experiments use the exponential map as a retraction. All manifolds were
implemented with or re-used existing implementations from the
\texttt{pymanopt}~\cite{JTownsend_NKoep_SWeichwald_2016} library in
\texttt{Python3}. The experiments were run on a computer with the
configuration: Debian GNU/Linux 13, Intel(R) Core(TM) i7-1265U 10 cores
4.8 GHz, 32 GB RAM. \rev{Running all methods} on the $100$ instances
of our problems took approximately three hours. \rev{For the manifold
	$\Mlinear$ (resp. $\Msphere$), computing an intrinsic PSS took $0.03m s$
	(resp. $0.1 ms$) on average while the computation of its projected counterpart
	required $0.2ms$ (resp. $1 ms$) on average. Although these values are
	implementation-dependent, they illustrate that the difference between the
	calculations is minimal, and justify further that we present profiles using
	the number of function evaluations as the budget.}

\subsection{Comparing projected and intrinsic variants}
\label{ssec:num:obs}

\begin{figure}[htb!]
	\centering
	\includegraphics[scale = 0.4]{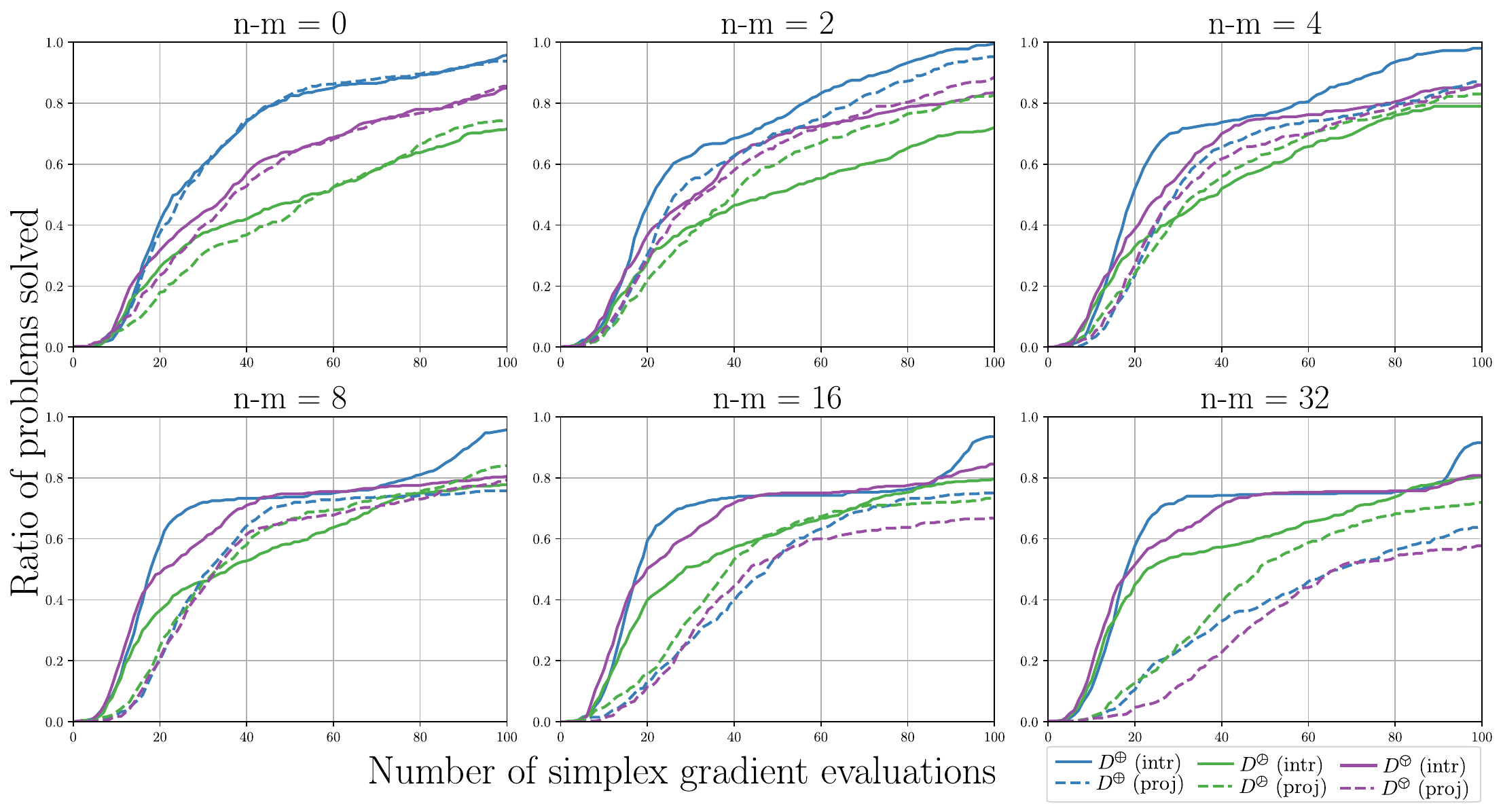}
	\caption{Data Profiles of PSSs $\pssa$, $\pssb$ and $\pssc$, intrinsic  variant \rrev{(plain)} or projected  variant \rrev{(dotted)}. $m = 4$. Each subplot corresponds \rev{to} a different codimension $n-m$. No rotation is applied. Tolerance: $\tau = 10^{-2}$.}
	\label{fig:num:codimsinfluence_norot}
\end{figure}

We first consider problems in our benchmark with $m=4$.
Figure~\ref{fig:num:codimsinfluence_norot} contains data profiles for all
possible codimensions. The intrinsic variant using $\pssa$ outperforms its
projected counterpart. \rev{The comparison is less pronounced for other
	choices of PSSs. Still, all intrinsic variants consistently outperform their
	projected counterparts as the codimension increases ($n-m=32$).}

\begin{figure}[htb!]
	\centering
	\includegraphics[scale = 0.4]{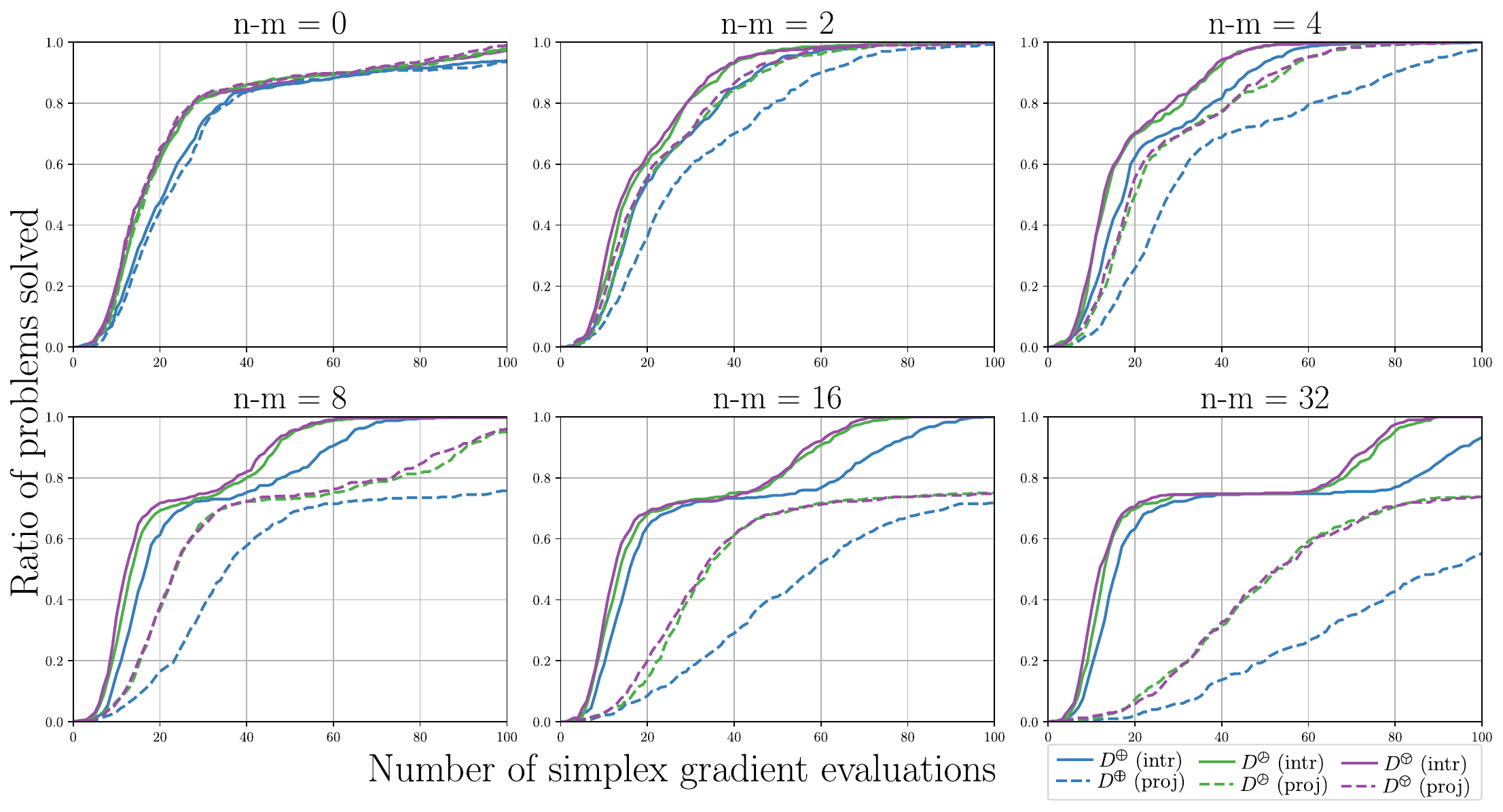}
	\caption{Data Profiles of PSSs $\pssa$, $\pssb$ and $\pssc$, intrinsic  variant \rrev{(plain)} or projected  variant \rrev{(dotted)}. $m = 4$. Each subplot corresponds to a different manifold codimension $n-m$. Rotation is applied. Tolerance: $\tau = 10^{-2}$.}
	\label{fig:num:codimsinfluence_rot}
\end{figure}

To determine whether using fixed PSSs aligned with coordinates directions (or a
specific basis) biases our results, we randomly rotate the polling sets at
every iteration of our methods. For the projected PSSs, this amounts to using
$\pssa(B)$, $\pssb(B)$, $\pssc(B)$, where $B$ is a random $n$ by $n$
orthogonal matrix, in lieu of $\pssa$, $\pssb$, $\pssc$. For the intrinsic
PSSs, it corresponds to using a random family of tangent space bases in lieu
of canonical ones. The results are given in
Figure~\ref{fig:num:codimsinfluence_rot}, and show a more consistent picture
than that of Figure~\ref{fig:num:codimsinfluence_norot}. Although the curves
are quite close for $n-m=0$ (admittedly a degenerate case), they show a
clearer advantage for intrinsic variants as the codimension \rev{increases}. We
observe in passing that using random rotations \rev{\rrev{seems} to worsen the performance of $\pssa$-based variants while improving that of other variants. This phenomenon has been previously observed in the
	Euclidean setting~\cite{CAudet_AIanni_SLeDigabel_CTribes_2014,
		SGratton_CWRoyer_LNVicente_2016}, while standard benchmarks have been shown
	to be biased towards coordinate directions~\cite{WHare_CWRoyer_2023}.}

\begin{figure}[htb!]
	\centering
	\includegraphics[scale = 0.4]{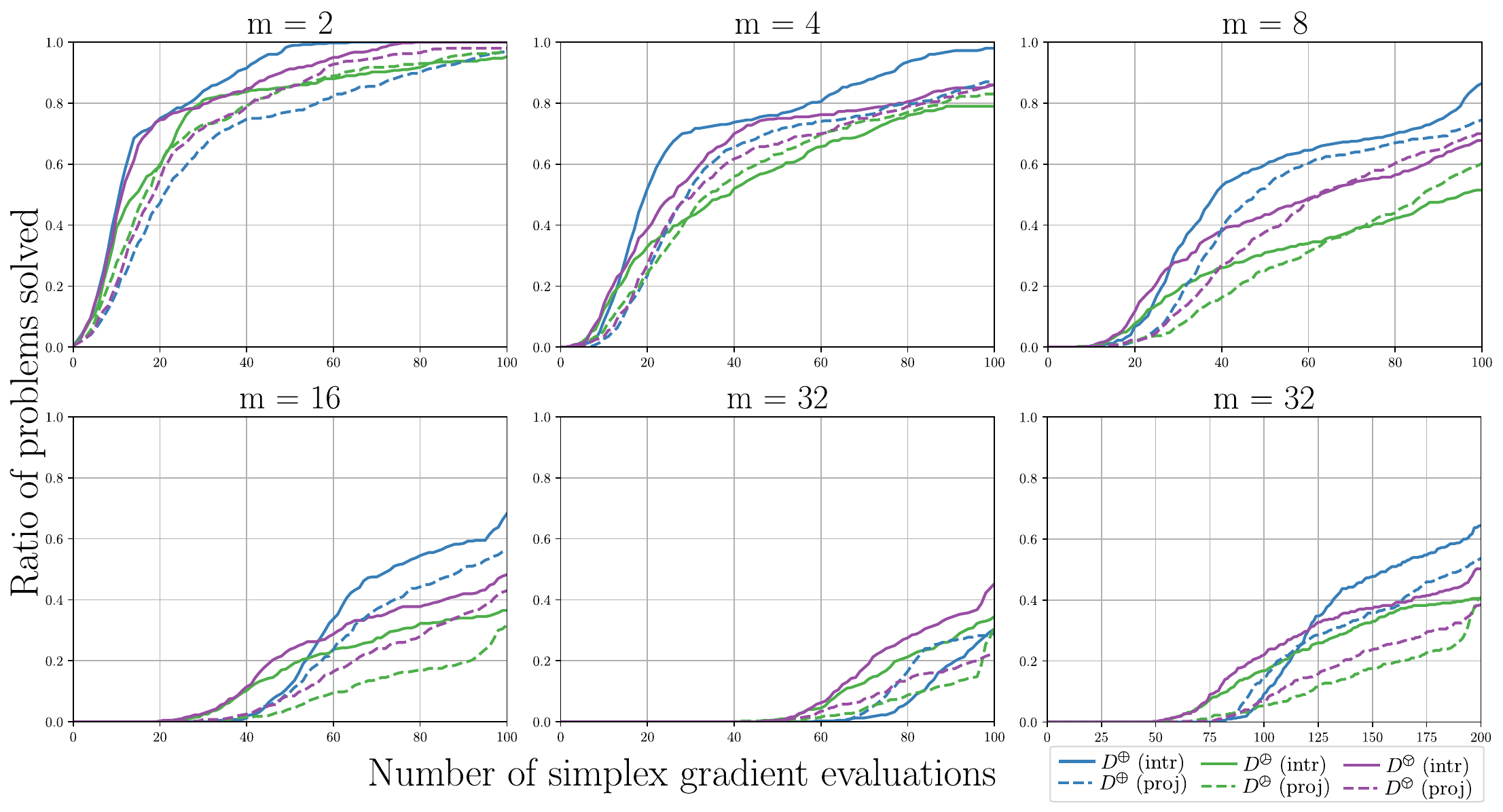}
	\caption{Data Profiles of PSSs $\pssa$, $\pssb$ and $\pssc$, intrinsic  variant \rrev{(plain)} or projected  variant \rrev{(dotted)}. $n-m = 4$. Each subplot corresponds to a different manifold dimension $m$. No rotation is applied. Tolerance: $\tau = 10^{-2}$. \rev{Lower right: budget of 200 simplex gradients.}}
	\label{fig:num:mdimsinfluence_norot}
\end{figure}

\begin{figure}[htb!]
	\centering
	\includegraphics[scale = 0.4]{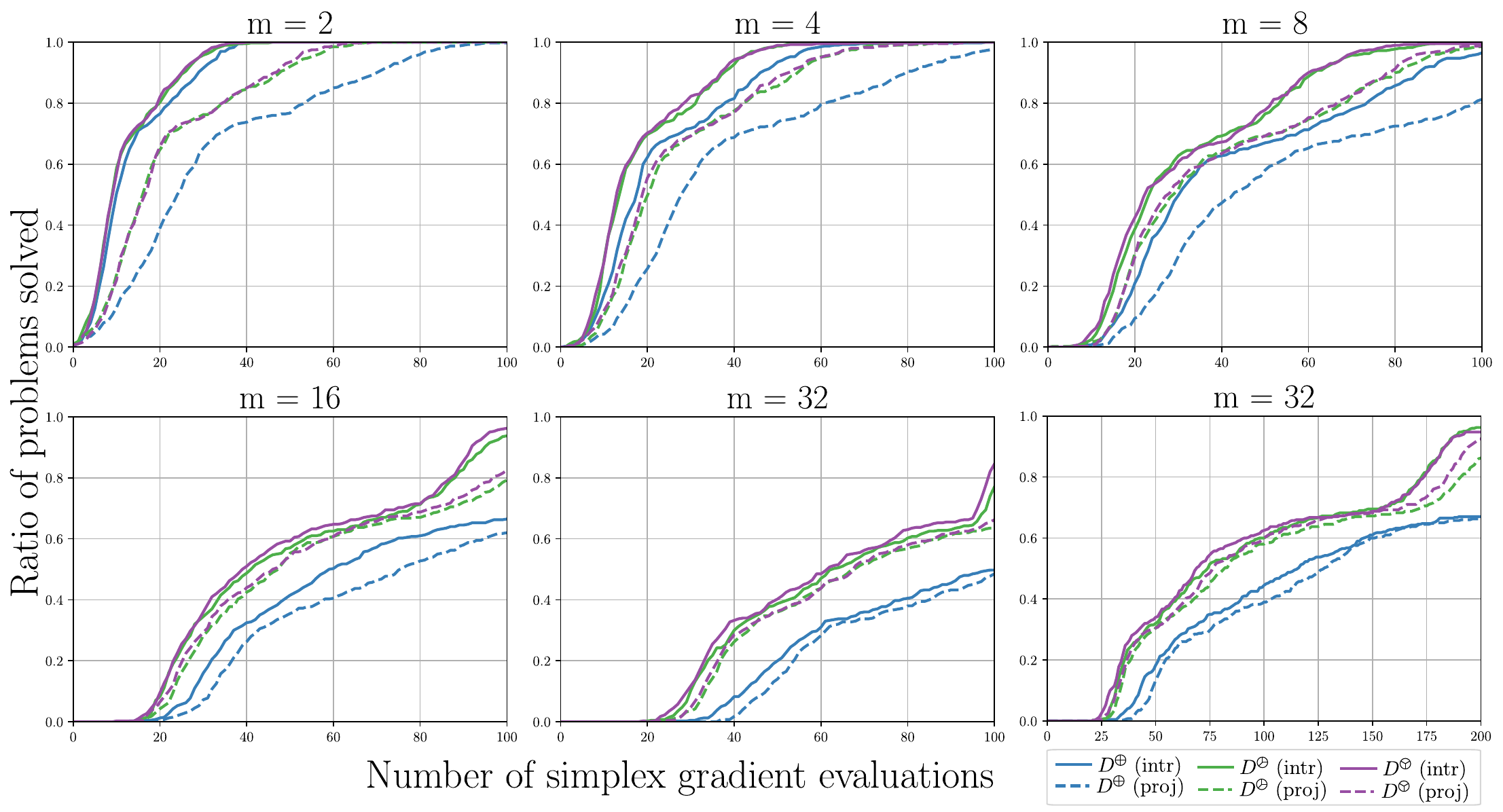}
	\caption{Data Profiles of PSSs $\pssa$, $\pssb$ and $\pssc$, intrinsic  variant \rrev{(plain)} or projected  variant \rrev{(dotted)}. $n-m = 4$. Each subplot corresponds to a different manifold dimension $m$. Rotation is applied. Tolerance: $\tau = 10^{-2}$. \rev{Lower right: budget of 200 simplex gradients.}}
	\label{fig:num:mdimsinfluence_rot}
\end{figure}

We now fix the codimension $n-m=4$. Figure~\ref{fig:num:mdimsinfluence_norot}
contains data profiles for all possible manifold dimensions $m$. Intrinsic methods
still exhibit better performance than their projected counterparts, but the picture
is less clear for intermediate dimensions $m = 4, 8$, and for the PSS variants using
$\pssa$. Adding random rotations to the polling directions again shows more
conclusive results (Figure~\ref{fig:num:mdimsinfluence_rot}). The gap between the
projected and intrinsic variants narrows as $m$ increases, yet using intrinsic
directions appears to be preferable.

We confirm these observations using Tables~\ref{table:pertable_norot}
and~\ref{table:perftable_rot}, that show the fraction of problems for which an
intrinsic method produces a better final function value than the corresponding
projected method. The tables further advocate using intrinsic variants when the
manifold codimension is large and the manifold dimension is small, especially
when the PSSs are built from $\pssa$ and/or when rotations are applied.

\begin{table}[htb!]
	\small
	\centering
	\caption{Performance tables on the benchmark problems for each PSS type. \rev{No rotation is applied.} Blue columns are associated to the data profiles in Figure \ref{fig:num:codimsinfluence_norot}\rev{. Red rows are associated to the data profiles in Figure \ref{fig:num:mdimsinfluence_norot}.}}
	\label{table:pertable_norot}
	\resizebox{\textwidth}{!}{%
		\subfloat[PSS  \protect $\pssa$]{%
			\hspace{.5cm}%
			\input{./rotation0_psstype1.tex}			\hspace{.5cm}%
		}
		\subfloat[PSS  \protect $\pssb$]{%
			\hspace{.5cm}%
			\input{./rotation0_psstype2.tex}			\hspace{.5cm}%
		}
		\subfloat[PSS  \protect $\pssc$]{%
			\hspace{.5cm}%
			\input{./rotation0_psstype3.tex}			\hspace{.5cm}%
		}%
	}%
\end{table}

\begin{table}[htb!]
	\small
	\centering
	\caption{Performance tables on the benchmark problems for each PSS type. \rev{Rotation is applied.} Blue columns are associated to the data profiles in Figure \ref{fig:num:codimsinfluence_rot}. Red rows are associated to the data profiles in Figure \ref{fig:num:mdimsinfluence_rot}.}
	\label{table:perftable_rot}
	\resizebox{\textwidth}{!}{%
		\subfloat[PSS \protect $\pssa$]{%
			\hspace{.5cm}%
			\input{./rotation1_psstype1.tex}			\hspace{.5cm}%
		}
		\subfloat[PSS  \protect $\pssb$]{%
			\hspace{.5cm}%
			\input{./rotation1_psstype2.tex}			\hspace{.5cm}%
		}
		\subfloat[PSS  \protect $\pssc$]{%
			\hspace{.5cm}%
			\input{./rotation1_psstype3.tex}			\hspace{.5cm}%
		}%
	}%
\end{table}



\rev{
	\subsection{Synchronization of rotations} \label{ssec:num:rotsync}

	Our previous experiments considered fairly simple manifolds in order to vary
	their dimension or codimension. We now turn to a more classical
	example of matrix manifold, \prev{namely} 
	\begin{equation*}
		\Mrotsync := \{(R_1, R_2) \in  \SO(5) \times \SO(5)\},
	\end{equation*}
	where $\SO(5) = \{R \in \R ^{5 \times 5}:R^TR = \rrev{I_5} ~\text{and}~ \det(R) = +1\}$
	is the special orthogonal group in $\R^{5 \times 5}$. This embedded manifold in
	$\R^{5 \times 5} \times \R^{5 \times 5}$ inherits its metric from the
	Frobenius inner product, and has dimension $m = 2\frac{5(5-1)}{2} = 20$ and
	codimension $n-m =50-20 = 30$. Given $(\hat R_1, \hat R_2) \in \Mrotsync$, the
	tangent space of $\Mrotsync$ at $(\hat R_1, \hat R_2)$
	is~\cite[Equation (7.37)]{NBoumal_2023}:
	\begin{equation}
		\label{eq:rot_sync_tangentspace}
		\T_{(\hat R_1,\hat R_2)}\Mrotsync
		=
		\hat R_1 \Skew(5) \times \hat R_2 \Skew(5),
	\end{equation}
	where $\Skew(5)$ is the set of skew-symmetric matrices in $\R^{5 \times 5}$.

	We consider an instance of synchronization of rotations, a problem aiming at
	recovering the relative orientations of multiple measurements of an
	object~\cite[Section 2.5]{NBoumal_2023}. Given a ground truth $(R_1,R_2)$
	(drawn randomly in $\Mrotsync$), we build a noisy version
	$H_{12} := R_1 R_2^T + \frac{10^{-6}}{\sqrt{5}} G$ of $R_1 R_2^T$, where
	$G \sim \mathcal{N}( 0,I_5)$. We then define the problem
	\begin{equation}
		\label{eq:rotsyncproblem}
		\min_{(\hat R_1,\hat R_2) \in \Mrotsync} \| \hat R_1 - H_{12} \hat R_2\|^2.
	\end{equation}

	We build projected PSSs of type $\pssa, \pssb$ and $\pssc$ as described in
	Section~\ref{ssec:projected_PSS} using the canonical basis of $\R^{50}$, since
	this space is isometric with $\R ^{5 \times 5}\times  \R ^{5 \times 5}$.
	Projected PSSs are computed by adapting the projection
	formula~\cite[Equation (7.38)]{NBoumal_2023} to the product
	manifold~\eqref{eq:rot_sync_tangentspace}. We build intrinsic PSSs using the
	canonical \rev{basis} of $\Skew(5)$ and the isometric isomorphism
	induced by~\eqref{eq:rot_sync_tangentspace} to construct bases for
	$\T_{(\hat R_ 1, \hat R_2)}\Mrotsync$. As in the previous section, we
	allow random rotations to be applied to these PSSs.
	Figure~\ref{fig:num:rotsync} shows data profiles obtained from solving $100$
	instances of problem~\eqref{eq:rotsyncproblem} using a budget of $500$
	Riemannian simplex gradients, \textit{i.e.} $500(20+1)$ function evaluations,
	random initial points in $\Mrotsync$, and tolerances
	$\tau \in \{0.01,0.05,0.1\}$. The results still favour intrinsic PSSs over
	their projected counterparts. Moreover, similarly to the previous figures, we
	observe that $\pssa$ performs better than $\pssb$ or $\pssc$ in the absence of
	rotations, but worse when rotations are applied.

	\begin{figure}[htb!]
		\centering
		\includegraphics[scale = 0.4]{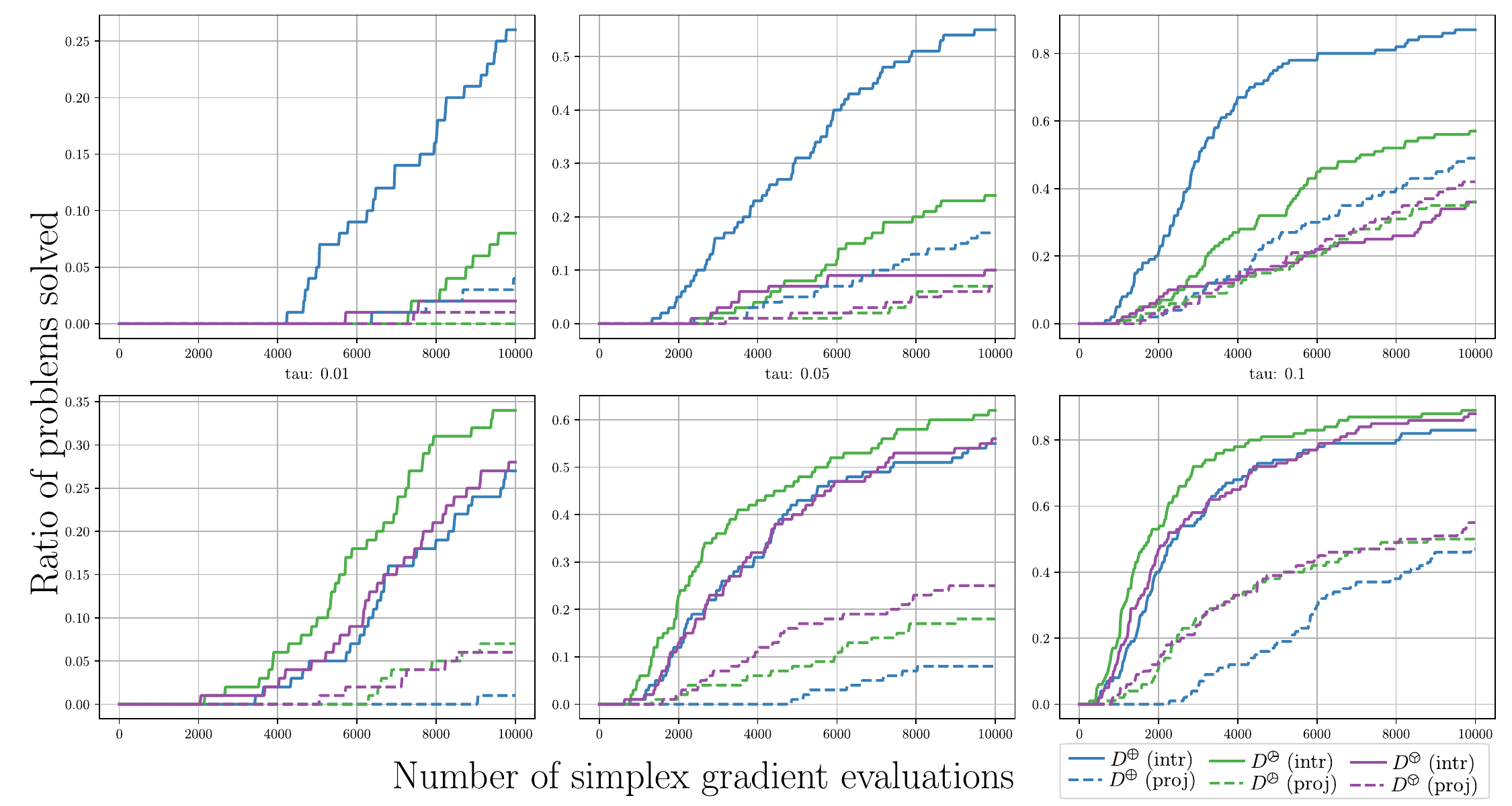}
		\caption{Data Profiles of PSSs $\pssa$, $\pssb$ and $\pssc$, intrinsic variant \rrev{(plain)} or projected  variant \rrev{(dotted)} for 100 instances of Problem \eqref{eq:rotsyncproblem}. Dimension $m=20$, codimension $n-m = 30$. First row: rotation is applied. Second row: rotation is not applied. Tolerance: $0.01$, $0.05$, $0.1$ (left to right).}
		\label{fig:num:rotsync}
	\end{figure}

}

\section{Conclusion}
\label{sec:conc}

In this paper, we derived a complexity analysis for a class of Riemannian
direct-search techniques. To this end, we generalized the Euclidean notion of
cosine measure to tangent spaces. Our complexity bounds depend on the cosine
measure and on the cardinality of polling sets used throughout the algorithm.
Our theory suggests that using intrinsic directions can be beneficial in
terms of both quantities, and our experiments suggest that the intrinsic
methods are generally more efficient than their projected counterparts, with
a stronger effect when the problem codimension is large with respect to the
manifold dimension.

Our analytical comparison of intrinsic and projected PSSs focused on the
hypersphere and the specific choice of PSS $\pssa$. Extending to other
manifolds, as well as other PSSs, is a natural continuation of our study.
Investigating probabilistic variants of direct-search, that have exhibited
better scalability than deterministic ones, also represents a future research
direction.

\paragraph{Funding} Funding for B. Cavarretta was provided
by \rev{Agence Nationale de la Recherche through the PEPR Compromis program
	\rev{22-PECY-0011} (France 2030)}. Funding for F. Goyens was provided by the
Fonds de la Recherche Scientifique - FNRS under Grant T.0001.23.
Funding for C. W. Royer was partially provided by Agence Nationale de la
Recherche through program ANR-23-IACL-0008 (PR[AI]RIE-PSAI), and by
CNRS IAE grant BONUS.

\bibliographystyle{plain}
\bibliography{dsmanifolds-revised-2}

\end{document}

%% file: rotation0_psstype1.tex
\begin{tabular}{l|ccccc}
\toprule
mdim & 2 & \textcolor{blue!70!black}{4} & 8 & 16 & 32 \\
codim &  &  &  &  &  \\
\midrule
0 & 0.52 & \textcolor{blue!70!black}{0.56} & 0.57 & 0.54 & 0.54 \\
2 & 0.94 & \textcolor{blue!70!black}{0.81} & 0.65 & 0.56 & 0.54 \\
\textcolor{red!70!black}{4} & \textcolor{red!70!black}{0.98} & \textcolor{purple!70!black}{0.95} & \textcolor{red!70!black}{0.76} & \textcolor{red!70!black}{0.62} & \textcolor{red!70!black}{0.52} \\
8 & 0.97 & \textcolor{blue!70!black}{0.97} & 0.87 & 0.70 & 0.52 \\
16 & 0.97 & \textcolor{blue!70!black}{0.97} & 0.95 & 0.85 & 0.64 \\
32 & 0.98 & \textcolor{blue!70!black}{0.97} & 0.97 & 0.92 & 0.77 \\
\bottomrule
\end{tabular}

%% file: rotation0_psstype2.tex
\begin{tabular}{l|ccccc}
\toprule
mdim & 2 & \textcolor{blue!70!black}{4} & 8 & 16 & 32 \\
codim &  &  &  &  &  \\
\midrule
0 & 0.58 & \textcolor{blue!70!black}{0.63} & 0.59 & 0.55 & 0.52 \\
2 & 0.53 & \textcolor{blue!70!black}{0.44} & 0.48 & 0.47 & 0.54 \\
\textcolor{red!70!black}{4} & \textcolor{red!70!black}{0.58} & \textcolor{purple!70!black}{0.44} & \textcolor{red!70!black}{0.46} & \textcolor{red!70!black}{0.47} & \textcolor{red!70!black}{0.48} \\
8 & 0.72 & \textcolor{blue!70!black}{0.46} & 0.41 & 0.46 & 0.50 \\
16 & 0.94 & \textcolor{blue!70!black}{0.66} & 0.49 & 0.50 & 0.55 \\
32 & 0.97 & \textcolor{blue!70!black}{0.76} & 0.66 & 0.55 & 0.56 \\
\bottomrule
\end{tabular}

%% file: rotation0_psstype3.tex
\begin{tabular}{l|ccccc}
\toprule
mdim & 2 & \textcolor{blue!70!black}{4} & 8 & 16 & 32 \\
codim &  &  &  &  &  \\
\midrule
0 & 0.57 & \textcolor{blue!70!black}{0.65} & 0.62 & 0.63 & 0.57 \\
2 & 0.65 & \textcolor{blue!70!black}{0.50} & 0.54 & 0.54 & 0.60 \\
\textcolor{red!70!black}{4} & \textcolor{red!70!black}{0.84} & \textcolor{purple!70!black}{0.50} & \textcolor{red!70!black}{0.54} & \textcolor{red!70!black}{0.54} & \textcolor{red!70!black}{0.58} \\
8 & 0.98 & \textcolor{blue!70!black}{0.61} & 0.52 & 0.53 & 0.63 \\
16 & 1.00 & \textcolor{blue!70!black}{0.82} & 0.64 & 0.62 & 0.64 \\
32 & 1.00 & \textcolor{blue!70!black}{0.95} & 0.74 & 0.66 & 0.66 \\
\bottomrule
\end{tabular}

%% file: rotation1_psstype1.tex
\begin{tabular}{l|ccccc}
\toprule
mdim & 2 & \textcolor{blue!70!black}{4} & 8 & 16 & 32 \\
codim &  &  &  &  &  \\
\midrule
0 & 0.54 & \textcolor{blue!70!black}{0.62} & 0.54 & 0.55 & 0.49 \\
2 & 1.00 & \textcolor{blue!70!black}{0.94} & 0.82 & 0.71 & 0.58 \\
\textcolor{red!70!black}{4} & \textcolor{red!70!black}{1.00} & \textcolor{purple!70!black}{0.99} & \textcolor{red!70!black}{0.89} & \textcolor{red!70!black}{0.77} & \textcolor{red!70!black}{0.70} \\
8 & 1.00 & \textcolor{blue!70!black}{1.00} & 0.97 & 0.87 & 0.82 \\
16 & 1.00 & \textcolor{blue!70!black}{1.00} & 1.00 & 0.95 & 0.90 \\
32 & 1.00 & \textcolor{blue!70!black}{1.00} & 0.99 & 0.98 & 0.95 \\
\bottomrule
\end{tabular}

%% file: rotation1_psstype2.tex
\begin{tabular}{l|ccccc}
\toprule
mdim & 2 & \textcolor{blue!70!black}{4} & 8 & 16 & 32 \\
codim &  &  &  &  &  \\
\midrule
0 & 0.43 & \textcolor{blue!70!black}{0.56} & 0.52 & 0.50 & 0.47 \\
2 & 0.92 & \textcolor{blue!70!black}{0.91} & 0.78 & 0.68 & 0.61 \\
\textcolor{red!70!black}{4} & \textcolor{red!70!black}{0.99} & \textcolor{purple!70!black}{0.97} & \textcolor{red!70!black}{0.89} & \textcolor{red!70!black}{0.78} & \textcolor{red!70!black}{0.72} \\
8 & 1.00 & \textcolor{blue!70!black}{1.00} & 0.96 & 0.84 & 0.80 \\
16 & 1.00 & \textcolor{blue!70!black}{1.00} & 0.99 & 0.93 & 0.86 \\
32 & 1.00 & \textcolor{blue!70!black}{1.00} & 0.99 & 0.98 & 0.96 \\
\bottomrule
\end{tabular}

%% file: rotation1_psstype3.tex
\begin{tabular}{l|ccccc}
\toprule
mdim & 2 & \textcolor{blue!70!black}{4} & 8 & 16 & 32 \\
codim &  &  &  &  &  \\
\midrule
0 & 0.41 & \textcolor{blue!70!black}{0.56} & 0.52 & 0.52 & 0.51 \\
2 & 0.90 & \textcolor{blue!70!black}{0.91} & 0.83 & 0.71 & 0.61 \\
\textcolor{red!70!black}{4} & \textcolor{red!70!black}{1.00} & \textcolor{purple!70!black}{0.97} & \textcolor{red!70!black}{0.88} & \textcolor{red!70!black}{0.80} & \textcolor{red!70!black}{0.71} \\
8 & 1.00 & \textcolor{blue!70!black}{0.99} & 0.96 & 0.88 & 0.80 \\
16 & 1.00 & \textcolor{blue!70!black}{1.00} & 0.99 & 0.94 & 0.86 \\
32 & 1.00 & \textcolor{blue!70!black}{1.00} & 1.00 & 0.97 & 0.94 \\
\bottomrule
\end{tabular}

%% file: dsmanifolds-revised-2.bbl
\begin{thebibliography}{10}

\bibitem{PAAbsil_RMahony_RSepulchre_2008}
P.-A. Absil, R.~Mahony, and R.~Sepulchre.
\newblock {\em Optimization Algorithms on Matrix Manifolds}.
\newblock Princeton University Press, Princeton, NJ, 2008.

\bibitem{AAgarwal_NBoumal_BBullins_CCartis_2021}
A.~Agarwal, N.~Boumal, B.~Bullins, and C.~Cartis.
\newblock Adaptive regularization with cubics on manifolds.
\newblock {\em Math. Program.}, 188:85--134, 2021.

\bibitem{CAudet_2014}
C.~Audet.
\newblock A survey on direct search methods for blackbox optimization and their
  applications.
\newblock In P.~Pardalos and T.~Rassias, editors, {\em Mathematics without
  boundaries}, pages 31--56. Springer, New York, NY, 2014.

\bibitem{CAudet_WHare_2017}
C.~Audet and W.~Hare.
\newblock {\em Derivative-{F}ree and {B}lackbox {O}ptimization}.
\newblock Springer Series in Operations Research and Financial Engineering.
  Springer International Publishing, Cham, Switzerland, 2017.

\bibitem{CAudet_WHare_GJarryBolduc_2025}
C.~Audet, W.~Hare, and G.~{Jarry-Bolduc}.
\newblock The cosine measure relative to a subspace.
\newblock {\em Comput. Optim. Appl.}, 92:125--153, 2025.

\bibitem{CAudet_AIanni_SLeDigabel_CTribes_2014}
C.~Audet, A.~Ianni, S.~{Le Digabel}, and C.~Tribes.
\newblock Reducing the number of function evaluations in mesh adaptive direct
  search algorithms.
\newblock {\em SIAM J. Optim.}, 24:621--642, 2014.

\bibitem{NBoumal_2023}
N.~Boumal.
\newblock {\em An introduction to optimization on smooth manifolds}.
\newblock Cambridge University Press, Cambridge, United Kingdom, 2023.

\bibitem{NBoumal_PAAbsil_CCartis_2019}
N.~Boumal, P.-A. {Absil}, and C.~Cartis.
\newblock Global rates of convergence for nonconvex optimization on manifolds.
\newblock {\em IMA J. Numer. Anal.}, 39:1--33, 2019.

\bibitem{ARConn_KScheinberg_LNVicente_2009}
A.~R. Conn, K.~Scheinberg, and L.N. Vicente.
\newblock {\em Introduction to derivative-free optimization}.
\newblock SIAM, Philadelphia, PA, 2009.

\bibitem{CDavis_1954}
C.~Davis.
\newblock Theory of positive linear dependence.
\newblock {\em Amer. J. Math.}, 76:733--746, 1954.

\bibitem{MDodangeh_LNVicente_ZZhang_2016}
M.~Dodangeh, L.~N. {Vicente}, and Z.~Zhang.
\newblock On the optimal order of worst case complexity of direct search.
\newblock {\em Optim. Lett.}, 10:699--708, 2016.

\bibitem{DWDreisigmeyer_2006}
D.~W. Dreisigmeyer.
\newblock Equality constraints, {R}iemannian manifolds and direct-search
  methods.
\newblock Technical Report LA-UR-06-7406, Los Alamos National Laboratory, 2006.

\bibitem{DWDreisigmeyer_2007a}
D.~W. {Dreisigmeyer}.
\newblock Direct search algorithms over {R}iemannian manifolds.
\newblock Technical Report LA-UR-06-7416, Los Alamos National Laboratory, 2007.

\bibitem{KJDzahini_FRinaldi_CWRoyer_DZeffiro_2025}
K.~J. {Dzahini}, F.~Rinaldi, C.~W. {Royer}, and D.~Zeffiro.
\newblock Direct-search methods in the year 2025: {T}heoretical guarantees and
  algorithmic paradigms.
\newblock {\em Euro. J. Comput. Optim.}, 13:100110, 2025.

\bibitem{SGratton_CWRoyer_LNVicente_2016}
S.~Gratton, C.~W. {Royer}, and L.~N. {Vicente}.
\newblock A second-order globally convergent direct-search method and its
  worst-case complexity.
\newblock {\em Optimization}, 65:1105--1128, 2016.

\bibitem{SGratton_CWRoyer_LNVicente_ZZhang_2019}
S.~Gratton, C.~W. {Royer}, L.~N. {Vicente}, and Z.~Zhang.
\newblock Direct search based on probabilistic feasible descent for bound and
  linearly constrained problems.
\newblock {\em Comput. Optim. Appl.}, 72:525--559, 2019.

\bibitem{WHare_GJarryBolduc_2020}
W.~Hare and G.~{Jarry-Bolduc}.
\newblock A deterministic algorithm to compute the cosine measure of a finite
  positive spanning set.
\newblock {\em Optim. Lett.}, 14:1305--1316, 2020.

\bibitem{WHare_GJarryBolduc_SKerleau_CWRoyer_2024}
W.~{Hare}, G.~{Jarry-Bolduc}, S.~Kerleau, and C.~W. {Royer}.
\newblock Using orthogonally structured positive bases for constructing
  positive $k$-spanning sets with cosine measure guarantees.
\newblock {\em Linear Algebra Appl.}, 680:183--207, 2024.

\bibitem{WHare_GJarryBolduc_CPlaniden_2023}
W.~Hare, G.~Jarry-Bolduc, and C.~Planiden.
\newblock Nicely structured positive bases with maximal cosine measure.
\newblock {\em Optim. Lett.}, 17:1495--1515, 2023.

\bibitem{WHare_CWRoyer_2023}
W.~Hare and C.~W. {Royer}.
\newblock Detecting negative eigenvalues of exact and approximate {H}essian
  matrices in optimization.
\newblock {\em Optim. Lett.}, 17:1739--1756, 2023.

\bibitem{WHare_SSun_2025}
W.~Hare and S.~Sun.
\newblock On the computation of the cosine measure in high dimensions.
\newblock arXiv:2506.19723, 2025.

\bibitem{GJarryBolduc_2023}
G.~Jarry-Bolduc.
\newblock {\em Numerical analysis for derivative-free optimization}.
\newblock PhD thesis, University of British Columbia, Vancouver, BC, Canada,
  2023.

\bibitem{TGKolda_RMLewis_VTorczon_2003}
T.~G. {Kolda}, R.~M. {Lewis}, and V.~Torczon.
\newblock Optimization by direct search: {N}ew perspectives on some classical
  and modern methods.
\newblock {\em SIAM Rev.}, 45:385--482, 2003.

\bibitem{VKungurtsev_FRinaldi_DZeffiro_2024}
V.~Kungurtsev, F.~Rinaldi, and D.~Zeffiro.
\newblock Retraction-based direct search methods for derivative-free
  {R}iemannian optimization.
\newblock {\em J. Optim. Theory Appl.}, 203:1710--1735, 2024.

\bibitem{SLeDigabel_SMWild_2024}
S.~{Le Digabel} and S.~M. {Wild}.
\newblock A taxonomy of constraints in black-box simulation-based optimization.
\newblock {\em Optim. Eng.}, 25:1125--1143, 2024.

\bibitem{JMLee_2018}
J.~M. {Lee}.
\newblock {\em Introduction to Riemannian manifolds}.
\newblock Grad. Texts in Math. Springer, Cham, Switzerland, second edition,
  2018.

\bibitem{LLiu_XZhang_2006}
L.~Liu and X.~Zhang.
\newblock Generalized pattern search methods for linearly equality constrained
  optimization problems.
\newblock {\em Appl. Math. Comput.}, 181:527--535, 2006.

\bibitem{JJMore_SMWild_2009}
J.~J. Mor\'{e} and S.~M. Wild.
\newblock Benchmarking derivative-free optimization algorithms.
\newblock {\em SIAM J. Optim.}, 20(1):172--191, 2009.

\bibitem{GNaevdal_2019}
G.~N{\ae}vdal.
\newblock Positive bases with maximal cosine measure.
\newblock {\em Optim. Lett.}, 13:1381--1388, 2019.

\bibitem{XPennec_2006}
X.~Pennec.
\newblock Intrinsic statistics on {R}iemannian manifolds: Basic tools for
  geometric measurements.
\newblock {\em J. Math. Imag. Vision}, 25(1):127--154, 2006.

\bibitem{RGRegis_2016}
R.~G. Regis.
\newblock On the properties of positive spanning sets and positive bases.
\newblock {\em Optim. Eng.}, 17(1):229--262, 2016.

\bibitem{RGRegis_2021}
R.~G. Regis.
\newblock On the properties of the cosine measure and the uniform angle
  subspace.
\newblock {\em Comput. Optim. Appl.}, 78(3):915--952, 2021.

\bibitem{RTRockafellar_1970}
R.~T. Rockafellar.
\newblock {\em Convex analysis}.
\newblock Princeton University Press, Princeton, NJ, 1970.

\bibitem{JTownsend_NKoep_SWeichwald_2016}
J.~Townsend, N.~Koep, and S.~Weichwald.
\newblock Pymanopt: {A} {P}ython toolbox for optimization on manifolds using
  automatic differentiation.
\newblock {\em J. Mach. Learn. Res.}, 17:1--5, 2016.

\bibitem{LNVicente_2013}
L.~N. {Vicente}.
\newblock Worst case complexity of direct search.
\newblock {\em Euro. J. Comput. Optim.}, 1:143--153, 2013.

\end{thebibliography}
